\documentclass[authoryear,preprint,11pt]{elsarticle}
\usepackage{amssymb}
\usepackage{eurosym}
\usepackage{color}
\usepackage{bm}
\usepackage{booktabs}
\usepackage{multirow}
\usepackage{longtable}
\usepackage{arydshln}
\usepackage{geometry}
\usepackage{enumitem}
\usepackage{comment}

\usepackage{subcaption}
\usepackage{tikz}
\usepackage{pgfplots}
\usetikzlibrary{patterns}
\usetikzlibrary{arrows}

\usepackage{amsthm,amsmath,eurosym}

\newtheorem{definition}{Definition}

\usepackage{lineno, blindtext}

\usepackage{etoolbox}
\makeatletter
\patchcmd{\ps@pprintTitle}{\footnotesize\itshape
Preprint submitted to \ifx\@journal\@empty Elsevier
\else\@journal\fi\hfill\today}{\relax}{}{}

\makeatother

\begin{document}

\begin{frontmatter}{}

\title{A multiple criteria nominal classification method in a web-based platform: Demonstration in a case of recruitment for the Portuguese Army}

\author[ceg-ist,inesc-id]{Ana Sara \textsc{Costa}\corref{cor1}}

\ead{anasaracosta@tecnico.ulisboa.pt}

\author[ceg-ist]{Jos{\'e} Rui \textsc{Figueira}}

\ead{figueira@tecnico.ulisboa.pt}

\author[inesc-id]{Jos{\'e} \textsc{Borbinha}}

\ead{jlb@tecnico.ulisboa.pt}

\address[ceg-ist]{CEG-IST, Instituto Superior T{\'e}cnico, Universidade de Lisboa, Portugal}

\address[inesc-id]{INESC-ID, Instituto Superior T{\'e}cnico, Universidade de Lisboa, Portugal}

\cortext[cor1]{Corresponding author. CEG-IST, Instituto Superior T{\'e}nico, Universidade
de Lisboa, Av. Rovisco Pais 1, 1049-001, Lisboa, Portugal. Tel.: +351 218 417 729. Date: April 8, 2019}

\begin{abstract}
\noindent 
\textsc{Cat-SD (Cat}egorization by Similarity-Dissimilarity) is a recently developed method for handling nominal classification problems
in the context of Multiple Criteria Decision Aiding (MCDA). This paper describes the design and implementation of this method, as well as an application dealing with a recruitment process in the Special Forces of the Portuguese Army. In addition, it proposes interaction protocols to elicit the preference parameters of the method to facilitate the construction of a decision model when the analyst guides the decision maker. \textsc{Cat-SD} has been implemented in \textsc{DecSpace}, a user-friendly on-line platform for supporting decision aiding processes  using one or more MCDA methods. The study related to the Portuguese Army Special Forces recruitment presents and demonstrates how these protocols and a tool like \textsc{DecSpace}
can facilitate the process of applying the method in real-world scenarios.
\end{abstract}
\begin{keyword}
 Decision Support Systems \sep Multiple Criteria Decision Aiding \sep Classification problems \sep Web-based platform \sep Special Forces.
\end{keyword}

\end{frontmatter}{}

\section{Introduction}\label{sec: Introduction}
\noindent Multiple Criteria Decision Aiding (MCDA) is a sub-discipline of Operations Research focus on the development and implementation of methods, techniques and tools that explicitly consider multiple criteria to assist a Decision Maker (DM) for dealing with decision situations. Some MCDA methods can be quite complex and some decision situations might require the use of more than one method. Therefore, having a Decision Support System (DSS) computing these methods and make it possible to visualize the results can be of an extreme valuable help. In fact, due to the increasing interest in this area, there has been a growing number of methods and, consequently, the number of DSS tools supporting them has been steadily increasing too. For example, two reviews of that are \cite{weistroffer2016software}, which offers not only a comprehensive survey, but also provides insights for future developments of software for this purpose, and \cite{mustajoki2017comparison}, motivated by the support to environmental planning processes, provide not only an analysis of well-known MCDA software, but also identifies the relevant aspects for good practices and innovative features. In fact, these tools can be challenging for users who are not experts in MCDA. Most of the available software solutions have been developed to meet the requirements of the researchers within this area, usually making available only one or a few methods, making it very hard for others to use that software. For that reason, requirements related to usability and user friendly issues need to be taken into account in the development of this kind of software solution. This is in line with the trend of filling the gap between MCDA researchers and practitioners \citep[see, e.g.,][]{thokala2018stakeholder,voinov2016modelling}.

In addition, many MCDA methods also might require the support of an analyst specialized in that method (but not forcibly in the tools) to help the DM (who must be knowledgeable about the problem requiring a decision, but not forcibly about the MCDA method or the tool). Techniques and interaction protocols between the analyst and the DM need further research to properly support the application of the MCDA methods. This kind of protocols are useful to guide the dialogue between the analyst and the DM when gathering the information related to the opinion, judgments and preferences of the DM. Such information is then used to assign properly values to the parameters needed for building a preference decision model.

In this paper, we focus on \textsc{Cat-SD (Cat}egorization by Similarity-Dissimilarity)  \citep{Costaetal2018}, a recently proposed MCDA nominal classification method. The main objective of this paper is to provide tools for supporting the application of this method in real-world cases. With such a purpose in mind, we introduce several protocols to be adopted by the analyst and the DM for eliciting the preference parameters in a co-constructive way (i.e., through a collaborative dialogue between the analyst and the DM), and we present the implementation of the \textsc{Cat-SD} method in a novel web-based platform  $-$ \textsc{DecSpace} (Decision Space), which makes use of forefront technology, usability and user interaction techniques  \citep{DecSpace_CAPSI,DecSpace_Andre}. While designing this platform, the main requirements were that it should be user-friendly and accessible, since it not only intends to make it possible to efficiently use various MCDA methods, but also to facilitate the understanding of the methods by non-expert users. The application of \textsc{Cat-SD} using \textsc{DecSpace} may help a DM to better understand the decision problem at hand.

The \textsc{Cat-SD} method has been developed to model the preferences of a DM and obtain assignment results when facing with multiple criteria nominal classification problem, while focusing on similarity and dissimilarity concepts. In the presence of a problem of this kind, a set of actions (or alternatives) has to be assigned to a set of categories previously defined in a nominal way (no preference order exists among them), considering a set of criteria and according to the preferences of the DM. For modeling purposes, the method needs data and preference parameters to be used as input for the computations, in order to obtain the categorization of the actions. A per-criterion similarity-dissimilarity function is used to model such concepts considering the differences of performance of two actions on a given criterion. Other preference parameters have to be defined, per category, namely reference actions, weights, possibly criteria interaction coefficients, and likeness thresholds. For eliciting such preference parameters, techniques are required, as proposed in  \citet{Costaetal2018}. In this paper, we propose interaction protocols to be used in a co-constructive way between the analyst and the DM to support the elicitation of such parameters.

\textsc{DecSpace} is a web-based platform that aims to provide an easy-to-use and intuitive approach to apply any kind of MCDA method, independently of its complexity. It is inspired in \textit{diviz} \citep{bigaret:hal-00926569, MeyerBigaret2012diviz}, in the sense that it permits to construct workflows that may use multiple MCDA methods and data, so the outputs of a certain method can be used as inputs to another one, allowing to design complete decision aiding processes. \textsc{DecSpace} is suitable to explore solutions for simple or complex decision situations and aims to be user-friendly to either non-expert users or MCDA expert users, with researching or teaching purposes, or even for consulting purposes.

This paper is organized as follows. Section \ref{sec:CAT-SD} contains a brief overview of the \textsc{Cat-SD} method. Section \ref{sec:CaseStudy} introduces a case study. Section \ref{sec:Preference-information} proposes generic interaction protocols for eliciting the preference parameters used in this method. Section \ref{sec:interaction} describes the interaction between the analyst and the DM, and the preference parameters defined to build the decision model. Section \ref{sec:DecSpace} is devoted to the design and application of \textsc{Cat-SD} in the \textsc{DecSpace} platform, and to demonstrate the use of the method by presenting the data and the model built with the intervention of the DM. Section \ref{sec:lessons} presents the main lessons learned from the application of the proposed protocols and the use of \textsc{DecSpace} while conducting the case study. Section \ref{sec:Conclusions} provides some concluding remarks, and proposes lines for future research and developments.

\section{Overview of the {\sc{Cat-SD}} method} \label{sec:CAT-SD}
\noindent This section is devoted to briefly introduce the \textsc{Cat-SD} method, designed to handle MCDA nominal classification problems. We present the main notation and the steps of the method. Related technical aspects, namely computations and the assignment procedure, are provided in Appendix. Details about the method can be found in \citet{Costaetal2018}.

\subsection{Basic notation} \label{subsec:notation}
\noindent Beforehand, several data sets and preference parameters have to be built in order to apply the method. The basic data of a multiple criteria nominal classification problem is presented in the following two sets: $A=\{a_1,\ldots,a_{i},\ldots\}$, which is the set of actions (not necessarily completely known \textit {a priori}); and $G=\{g_{1},\ldots,g_{j},\ldots,g_{n}\}$, which is the a coherent family of criteria as defined according to \cite{royaide}. Let $g_j(a_i)$ denote the performance of action $a_i$ on criterion $g_j$. A performance table can be built. In addition, we consider a set of non-ordered categories, $C=\{C_{1}\ldots,C_{h},\ldots,C_{q},C_{q+1}\}$, also called nominal categories ($C_{q+1}$ is a dummy category that receives actions which are not assigned to the other categories).

The application of the method needs the construction of several preference parameters. For characterizing the categories, we consider $B=\{B_{1},\ldots,B_{h},\ldots,B_{q+1}\}$ as the set of all reference actions that allow to defined such categories (where $B_{q+1}=\emptyset$), and  $B_{h}=\{b_{h1},\ldots,b_{h\ell},\ldots,b_{h\vert B_{h}\vert}\}$ is the set of reference actions considered to define category $C_{h}$, for $h=1,...,q$.

In addition, we must consider the following preference parameters. 
\begin{itemize}[label={--}]
\item $k_{j}^{h}$ is the weight of criterion $g_{j}$ (its relative importance) with respect to category $C_{h}$, for $j=1,\ldots,n$ and $h=1,\ldots,q$;
\item $k_{j\ell}^{h}$ is the mutual-strengthening (or mutual-weakening) coefficient of the criteria pair $\{g_{j},g_{\ell}\}$, with $k_{j\ell}^{h}>0$ (or $k_{j\ell}^{h}<0$), for $h=1,\ldots,q$;
\item $k_{j|p}^{h}$ is the antagonistic coefficient for the ordered criteria pair $(g_{j},g_{p})$, with $k_{j|p}^{h}<0$, for $h=1,\ldots,q$;
\item $\lambda^{h}$ is the likeness threshold defined for category $C_{h}$, for $h=1,\ldots,q$.
\end{itemize}

\subsection{The flowchart of the method} \label{subsec:flowchart}
\noindent  The following flowchart contains three main blocks: one devoted to the input, other to the calculations, and the last is related to the output provided by the method (i.e., the assignment results). They can be briefly described as follows.

\begin{enumerate}
    \item \textit{Input.} There are basically two types of input:
    \begin{enumerate}
        \item \textit{Data}. The basic data is composed of the set of criteria $G$, the set of actions $A$, the performance table, and the set of categories $C$;
        \item \textit{Preference parameters}. These are in the second box of the input block, in the flowchart, and they are the set of reference actions $B$, the weights of criteria and the interaction coefficients per category. Per-criterion similarity-dissimilarity functions (SD functions) have to be constructed. In addition, a likeness threshold per category should also be defined. The notation for these preference parameters is provided in the previous subsection; 
    \end{enumerate}
    \item \textit{Calculations.} The calculation phase is composed of a sequence of steps as presented in Appendix, Steps 1-5;
    \item \textit{Output.} As output, the method provides the assignment of the actions to the non-ordered categories. An action can be assigned to one or several categories, or not assigned at all (i.e., it is assigned to the dummy category).
\end{enumerate}

\begin{figure}[htb!]
\centering 
\includegraphics{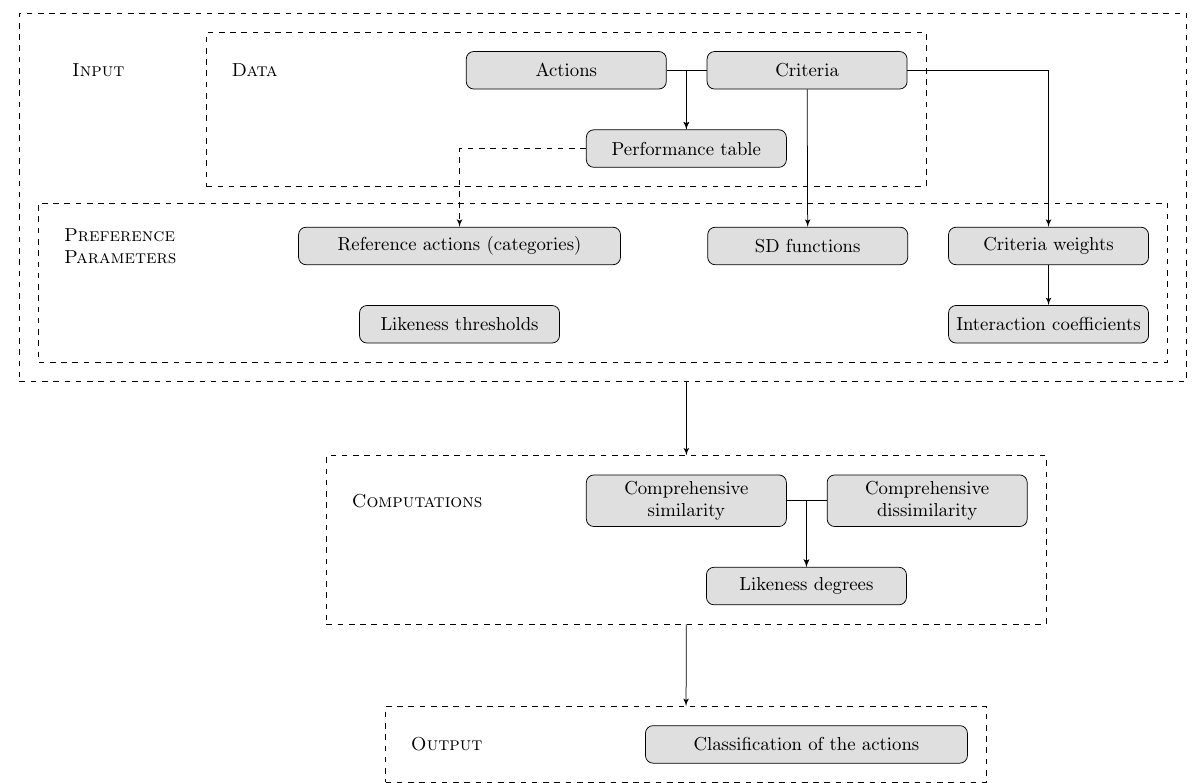} 
\caption{A flowchart of the \textsc{Cat-SD} method}
\label{Fig_flowchart}
\end{figure}

The next subsection presents a case study that  help us to show how to proceed to transform the preference information into parameters as proposed in Section \ref{sec:Preference-information}.

\section{Case study} \label{sec:CaseStudy}
\noindent This section introduces a case study, comprising its context and the required data to construct a \textsc{Cat-SD} decision model.

The Portuguese Armed Forces are a national institution in charge of guarantee the national independence and the territorial integrity of Portugal. The Portuguese Army, one of three branches (the largest branch), contributes to this mission as landing forces, in cooperation with the Portuguese Navy and Air Force. The Portuguese Army is also an institution as old as the national foundation, having been always present in great moments of the history of Portugal. 

Currently, the Portuguese Army essentially generates forces and capabilities, for projection and employment inside and outside the country's borders, to support people
and their institutions. In order to achieve its mission it needs adequate resources, being the human factor, around twelve thousand military
and civilians, its greatest valuable resource and the primary cause of its efficacy.

The soldier is the principal value of the military forces, and, for this reason, the military recruitment should be an essential and permanent concern. Today, the Portuguese military have been operated in demanding theater of operations, namely in Iraq, Afghanistan, and Central African Republic, in order to defend the State interests under the scope of cooperative security. In these operational environments, complex and presenting high risk, it is used Special Forces with advance training and highly qualified. The Special Forces that integrate the Portuguese Army are mainly divided in Commandos, Paratroopers and Special Operations. Included in the Special Operations are the snipers, being elite forces recognized and praised by their high performance.

Accordingly, the process of recruitment and selection of military candidates for Special Forces requires particular attention. Upon completion of the initial military training, a military intending to join the Portuguese Army Special Forces has to fulfill some prerequisites. Depending on the type of Special Forces, different profiles are required, based on a set of features. Those features are evaluated through ability tests related to physical fitness, psychological aspects, and health and medical issues. Interviews are also applied in order to assess personality traits.

In the presented context, we have conducted a study, in collaboration with the Portuguese Army through the \textit{Centro de Psicologia Aplicada do Ex{\'e}rcito (CPAE)}, which is the military center evolved in the process of military recruitment with the mission of studying, applying and supervising the activities of the Army in the area of psychology. The DM is an expert in the assessment and selection of the candidates to the Army. 

The case study aims at assigning military candidates into a small number of categories, which represent distinct types of Portuguese Army Special Forces. For modelling purposes, the \textsc{Cat-SD} method was applied. For that, with intervention of the DM, in a co-constructive way, a decision model has been built, including a set of criteria and a set of preference parameters reflecting the judgments of the DM. For confidentially reasons, some data used in this study are not real.

In general, a decision aiding process is a long process requiring several interaction sessions between analysts and DMs, to end-up with a decision classification model. However, in this case, conducting this study was immediately well-accepted by the DM, since it is in  line with the current Portuguese Army recruitment process. Once the model is built, it can be used in the future and the DM is more open to have long interaction sessions when needed.

In the first meeting, the analyst and the DM formulated the classification problem, representing the decision situation at hand. The DM refereed that the profile of a sniper is very specific and, for that reason, a category should be defined for the sniper profile. Thus, in this study, four (nominal) categories have been considered (five categories when the dummy category, $C_{5}$, is counted), as described below.

\begin{itemize}[label={--}]
\item \textit{Commandos} ($C_{1}$): light unshielded combat forces intended for conventional operations of eminently offensive nature, taking advantage of the surprise, speed, violence and precision of the attack. They have great technical and tactical abilities, and high state of readiness;
\item \textit{Paratroopers} ($C_{2}$): military parachutists (light infantry
force) who have the capability of insertion in the theater of operations through parachute jumping, usually in surprise attacks. They are characterized
by a capability of a great speed in action and flexibility, presenting a high state of readiness;
\item \textit{Special Operations} ($C_{3}$): military that uses unconventional methods and resources to operate in diverse environments. They are generally adaptable and self-reliant;
\item \textit{Snipers} ($C_{4}$): military who shoot targets from long distances, using high-precision rifles and high-magnification optics, while maintaining effective visual contact with enemy targets from concealed positions. They typically present character traits of patience, attention to detail, perseverance, and physical endurance, and an attitude of determination;
\item \textit{Unsuitable candidates} ($C_{5}$): military that are not suitable to any previous category according to the constructed decision model.
\end{itemize}

It should be recalled that the application of \textsc{Cat-SD} provides at least a category for which a candidate should be assigned, meaning that a candidate may be suitable for more than one special force category according to the model. Based on the assignment results, a final decision is taken by the DM, and then each candidate assigned to a given category will be subject to a similar training course.

The following are the criteria considered relevant by the DM for the assignment of the military candidates into the Special Forces. Let us remark that all the criteria are to be maximized.

\begin{itemize}[label={--}]
\item \textit{Physical fitness} ($g_{1}$) - $PF$: It consists of the current level of fitness assessed by the performance of the candidate on some exercises, such as wall transposition, swimming, pushing-up, curling-up, ditch jumping and running. The result of each test is converted to a corresponding value within the range $[0,20.00]$, with exception of the result of the decision tests (e.g., wall transposition), which can be \textit{yes} or \textit{no}. For a candidate, the average of her/his numerical results of the tests is calculated, within the range $[0,20.00]$, while the result in the decision tests is \textit{yes}. Only candidates with an average score greater or equal to $10.00$ are considered. Thus, the performances of this criterion are assessed in a quantitative scale within the range $[10.00,20.00]$;
\item \textit{Intelligence} ($g_{2}$) - $Intel$: It includes intellectual aptitudes of the candidate, while evaluating the competency in logical reasoning, the capacity to understand, plan and solve problems. The performances of this criterion are assessed in a qualitative scale with five levels: insufficient ($1$), weak ($2$), satisfactory ($3$), good ($4$), excellent ($5$);
\item \textit{Numerical reasoning} ($g_{3}$) - $NR$: It consists of the candidate's ability to interpret, analyze and draw logical conclusions from numerical data and make reasoned decisions when solving problems. The performances of this criterion are assessed in a quantitative scale, corresponding to percentile scores, within the range $[30,99]$;
\item \textit{Spatial ability} ($g_{4}$) - $SA$:  It consists of the candidate's capacity of the spatial perception, including the ability to understand the spatial relations among objects and to mentally visualize and manipulate them. The performances of this criterion are assessed in a quantitative scale, through the application of a perception-cognitive test, resulting in a percentile score, within the range $[30,99]$;
\item \textit{Mechanical reasoning} ($g_{5}$) - $MechR$: It consists of the practical knowledge of mechanics and physics, that is, the candidate's capability to understand and apply the concepts and principles of mechanics within a variety of situations. The performances of this criterion are assessed in a quantitative scale, through the application of a perception-cognitive test, resulting in a percentile score, within the range $[30,99]$;
\item \textit{Velocity perception} ($g_{6}$) - $VP$: It is related to the processing speed of the candidate, a cognitive ability to do a mental task, i.e., the time to capture and react to the information received (the stimulus can be whether visual, auditory or movement). The performances of this criterion are assessed in a quantitative scale, through the application of a perception-cognitive test, resulting in a percentile score, within the range $[30,99]$;
\item \textit{Psychomotor ability} ($g_{7}$) - $PmA$: It consists of the physical movement related to conscious cognitive processing, associated to the motor coordination and manual dexterity. The performances of this criterion are assessed in a qualitative scale with five levels: very low ($1$), low ($2$), medium ($3$), high ($4$), very high ($5$);
\item \textit{Personality} ($g_{8}$) - $Pers$: It includes some candidate's personality traits and adaptation abilities in the military environment, such as emotional stability, maturity, adaptability, resilience, teamwork skills and motivation. The performances of this criterion are assessed in a qualitative scale with five levels: very weak ($1$), weak ($2$), medium ($3$), good ($4$), very good ($5$);
\item \textit{Physical condition (medical)} ($g_{9}$) - $Med$: It consists of the assessment of medical aspects, such as the candidate's physical constitution, audition and vision. The performances of this criterion are assessed in a qualitative scale with five levels: clearly bellow average ($1$), bellow average ($2$), average ($3$), above average ($4$), clearly above average ($5$).
\end{itemize}

A set of candidates to become Special Forces soldiers have been assessed according to the nine predefined criteria. Table \ref{Table:Performance} displays the performances on all criteria of twenty candidates (the data are not real).

\begin{table}[!htb]
\caption{Performance of the candidates to Special Forces on each criterion}
\label{Table:Performance}\smallskip{}
\centering \begin{small} \resizebox{0.75\textwidth}{!}{ %
\begin{tabular}{llllllllll}
\hline 
Candidate & \textbf{$PF$} & \textbf{$Intel$} & \textbf{$NR$} & \textbf{$SA$} & \textbf{$MechR$} & \textbf{$VP$} & \textbf{$PmA$} & \textbf{$Pers$} & \textbf{$Med$}\tabularnewline
\hline 
$a_{1}$  & 17.25 & 4 & 65 & 75 & 70 & 75 & 4 & 4 & 4\tabularnewline
$a_{2}$  & 16.05 & 4 & 85 & 85 & 90 & 80 & 4 & 5 & 5\tabularnewline
$a_{3}$  & 14.91 & 4 & 60 & 75 & 85 & 55 & 4 & 5 & 4\tabularnewline
$a_{4}$  & 15.00 & 3 & 65 & 85 & 80 & 65 & 4 & 4 & 5\tabularnewline
$a_{5}$  & 13.73 & 4 & 75 & 96 & 75 & 70 & 4 & 4 & 3\tabularnewline
$a_{6}$  & 18.28 & 3 & 70 & 75 & 60 & 70 & 4 & 5 & 4\tabularnewline
$a_{7}$  & 12.83 & 5 & 80 & 60 & 75 & 85 & 4 & 4 & 3\tabularnewline
$a_{8}$  & 14.50 & 4 & 75 & 80 & 96 & 80 & 5 & 5 & 5\tabularnewline
$a_{9}$  & 15.75 & 4 & 55 & 65 & 75 & 97 & 5 & 5 & 5\tabularnewline
$a_{10}$  & 15.86 & 4 & 90 & 80 & 75 & 80 & 2 & 5 & 4\tabularnewline
$a_{11}$  & 19.12 & 3 & 50 & 75 & 65 & 75 & 4 & 4 & 5\tabularnewline
$a_{12}$  & 14.35 & 2 & 80 & 85 & 85 & 70 & 4 & 3 & 4\tabularnewline
$a_{13}$  & 11.65 & 4 & 75 & 85 & 96 & 65 & 4 & 4 & 4\tabularnewline
$a_{14}$  & 16.00 & 5 & 80 & 55 & 65 & 75 & 3 & 4 & 4\tabularnewline
$a_{15}$  & 18.00 & 3 & 75 & 70 & 50 & 75 & 4 & 4 & 5\tabularnewline
$a_{16}$  & 17.22 & 4 & 60 & 70 & 75 & 85 & 3 & 4 & 5\tabularnewline
$a_{17}$  & 13.85 & 4 & 90 & 85 & 80 & 90 & 5 & 5 & 4\tabularnewline
$a_{18}$  & 15.10 & 3 & 70 & 90 & 95 & 60 & 5 & 4 & 4\tabularnewline
$a_{19}$  & 12.45 & 5 & 80 & 65 & 70 & 70 & 4 & 5 & 4\tabularnewline
$a_{20}$  & 14.32 & 4 & 85 & 80 & 85 & 75 & 5 & 5 & 5\tabularnewline
\hline 
\end{tabular}} \end{small} 
\end{table}

\section{Assessing the DM preferences and the modeling of the parameters} \label{sec:Preference-information}

\noindent As stated in \citet{Costaetal2018}, different approaches can be adopted to assess the preference parameters of our model, starting from a pure learning from data based approach with ``training''  examples provided by the DM to a full co-constructive approach. In this work, we focus on the latter, in which the analyst and the DM interact and cooperate to assess the preferences and model them in view to construct a preference model.  This section is devoted to present interaction protocols, in order to assign appropriate values to the parameters according to the preferences of the DM. This kind of protocols are based on a structured dialogue between the analyst (the questioner) and the DM (the questionee) involving introductory and easy questions \citep[Chapter 11]{roy1996multicriteria}.

\subsection{The reference actions}
\noindent Each category is defined by a non-empty set of reference actions. A reference action is a representative action  of a given category, i.e, a prototype or a typical element. In order to build the set of reference actions, the analyst can start by asking the DM to choose a category. If the DM is able o identify at least a representative action (e.g., from past decisions, guidelines, holistic judgments, etc.) of such a category, those actions can be used to characterize the considered category. If this is not possible, at least a dummy action must be built. This can be done by considering some adequate performances levels in the criteria scales. The analyst and the DM should proceed in this way while considering individually each one of the predefined categories. In the end, each category is defined by at least one reference action (with exception of the dummy category).

\subsection{The per-criterion functions} \label{subsec:SD_functions}
\noindent A \textit{per-criterion similarity-dissimilarity function}, $f_j(\Delta_j(a,b))$, is used for each criterion, $g_j$, to model similarity-dissimilarity
judgments in the comparison of each ordered pair $(a,b)$ of actions, where $\Delta_j(a,b)$ is the performance difference between $g_j(a)$ and $g_j(b)$. (for more details, see \citealp{Costaetal2018}). A general definition of this function is also presented in \ref{sec:appendix}.
The per-criterion similarity-dissimilarity functions make use of  similarity-dissimilarity thresholds are defined such that $v\big(g_{j}(b)\big)\geqslant u_{j}\big(g_{j}(b)\big)\geqslant t_{j}\big(g_{j}(b)\big)\geqslant0$
and $v_{j}^{'}\big(g_{j}(b)\big)\geqslant u_{j}^{'}\big(g_{j}(b)\big)\geqslant t_{j}^{'}\big(g_{j}(b)\big)\geqslant0$. If the difference of performances, $diff\{g_{j}(a),g_{j}(b)\}$,
is within the range $]-t_{j}^{'}\big(g_{j}(b)\big),t_{j}\big(g_{j}(b)\big)[$,
then there is a positive contribution to the similarity on the criterion
$g_{j}$. If the difference of performances, $diff\{g_{j}(a),g_{j}(b)\}$,
is within the ranges $[-diff\{g_{j}^{\max},g_{j}^{min}\},-u_{j}^{'}\big(g_{j}(b)\big)[$
and $]u_{j}\big(g_{j}(b)\big),diff\{g_{j}^{\max},$ $g_{j}^{min}\}]$, then there is a negative contribution to the similarity: a negative dissimilarity when $a$ is strictly less than $b$, and a positive dissimilarity when $a$ is strictly greater than $b$.

The thresholds can be constant thresholds, i.e, when they are invariable along the scale of the criterion, which means that the same value is used to compare two actions and does not depend on $g_{j}(b)$, while variable thresholds vary along the range of the criterion scale.

In the next subsection we provide a more formal definition of all the types of thresholds mentioned above.

\subsubsection{Definition of the similarity and dissimilarity thresholds\label{subsec:Building-thresholds-for}}
\noindent The  thresholds can be
defined as in Definition \ref{def:thresholds} below.

\begin{definition}[Similarity and dissimilarity thresholds]\label{def:thresholds}
Three different kind of thresholds are defined as follows, for two actions $a$ and $b$, where $a$ is the action to be assessed and $b$ is the reference action.

\begin{itemize}
\item[{$i)$}] The similarity thresholds, $t_{j}\big(g_{j}(b)\big)$ and $t_{j}^{'}\big(g_{j}(b)\big)$, can be defined as follows:
\begin{itemize}
\item Consider $g_{j}(a) \geqslant g_{j}(b)$. The threshold $t_{j}\big(g_{j}(b)\big)$
is the largest performance difference that allows to consider that
action $a$ is similar to action $b$ according to criterion $g_{j}$; 
\item Consider $g_{j}(a) \leqslant g_{j}(b)$. The threshold $t_{j}^{'}\big(g_{j}(b)\big)$
is the largest performance difference that allows to consider that
action $a$ is similar to action $b$ according to criterion $g_{j}$;
\end{itemize}
\item[{$ii)$}] The dissimilarity thresholds,
$u_{j}\big(g_{j}(b)\big)$ and $u_{j}^{'}\big(g_{j}(b)\big)$, can
be defined as follows: 
\begin{itemize}
\item Consider $g_{j}(a) \geqslant g_{j}(b)$. The threshold $u_{j}\big(g_{j}(b)\big)$
is the smallest performance difference that allows to consider that
action $a$ is dissimilar to action $b$ according to criterion $g_{j}$;
\item Consider $g_{j}(a)\leqslant g_{j}(b)$. The threshold $u_{j}^{'}\big(g_{j}(b)\big)$
is the smallest performance difference that allows to consider that
action $a$ is dissimilar to action $b$ according to criterion $g_{j}$;
\end{itemize}
\item[{$ii)$}] The total dissimilarity thresholds, $v_{j}\big(g_{j}(b)\big)$ and
$v_{j}^{'}\big(g_{j}(b)\big)$, can be defined as follows: 
\begin{itemize}
\item Consider $g_{j}(a) \geqslant g_{j}(b)$. The threshold $v_{j}\big(g_{j}(b)\big)$
is the smallest performance difference that allows to consider that
action $a$ is totally dissimilar to action $b$ according to criterion
$g_{j}$;
\item Consider $g_{j}(a)\leqslant g_{j}(b)$. The threshold $v_{j}^{'}\big(g_{j}(b)\big)$
is the smallest performance difference that allows to consider that
action $a$ is totally dissimilar to action $b$ according to criterion
$g_{j}$.
\end{itemize}
\end{itemize}

\end{definition}

\subsubsection{Assessing similarity and dissimilarity thresholds\label{subsec:Building-thresholds-for}}
\noindent  The assignment of values to these thresholds can be carried out in a constructive way adopting a similar protocol usually used in the elicitation of veto thresholds in \textsc{Electre} methods \citep{roy2014discriminating}.
The analyst can start by using some levels from the criterion scale as reference values to help DM in the assessment of the  thresholds values. According to Roy et al. (2014), for discrete scales, if checking all levels of the scale the conclusion is that the value of the threshold is the same, then we are in presence of a constant threshold. Otherwise, the threshold has to be defined for each level of the criterion scale. They should proceed in a analogous way with all thresholds.

\paragraph{Gathering the preference information}

For continuous scales, through a co-construction interactive process between the analyst and the DM a way of determining the thresholds can be as follows (see also Figure \ref{Fig:SD_function}):

\begin{figure}[htb!]
\centering 
\includegraphics[scale=0.65]{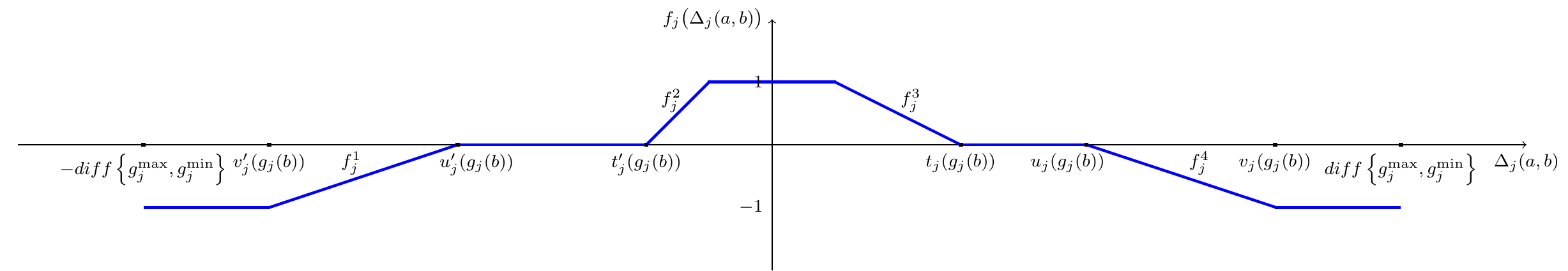} 
\caption{An example of a per-criterion similarity-dissimilarity function}
\label{Fig:SD_function}
\end{figure}

\begin{enumerate}
\item Consider a criterion, $g_{j}$, to be maximized. Choose a reference scale level in the lower part of the criterion scale. Let $g_{j}^{l}(b)$ denote such a level. In general, a \textquotedblleft good\textquotedblright{}
choice is such that the chosen level is in between $1/4$ to $1/3$ of the scale. \label{enu:Consider-a-criterion}
\item Consider the two performance scale levels, $g_{j}^{l}(b)$ and $g_{j}(a)$.
The latter should be closed enough to the reference level. In general, with $g_{j}\left(a\right)>g_{j}^{l}(b)$, but we can start with the equality.
The analyst starts by asking the DM whether, $g_{j}^{l}(b)$ and $g_{j}(a)$,
can be considered similar from her/his point of view: \label{enu:Consider-two-performance}
\begin{enumerate}
\item If the DM answers negatively, the analyst incrementally decreases
$g_{j}(a)$, and asks the same question until obtain a positive answer (if any).
This may occur only when $g_{j}^{l}(b)$ and $g_{j}(a)$ become equal.
Then, $g_{j}\left(a\right)$ is fixed. Let $g_{j}^{lt}\left(a\right)$ denote such a performance level. \label{enu:decrease_g(a)}
\item If the DM answers positively, the analyst incrementally increases
$g_{j}(a)$, and asks the same question until she/he gets a negative answer. Then, $g_{j}\left(a\right)$ is fixed as $g_{j}^{lt}\left(a\right)$.
\label{enu:increase_g(a)}
\end{enumerate}
The objective of this step is thus to identify, in a co-constructive way, i.e., through a dialogue between the analyst and the DM, how many (if any) performance levels of action $a$ must be removed from the performance, $g_{j}^{l}(b)$, to consider that both actions, $a$ and $b$, are similar, according to criterion $g_{j}$. The performance level $g_{j}^{lt}\left(a\right)$ is the lowest constructed performance level allowing to conclude about the similarity with respect to $g_{j}^{l}(b)$.
Any other performance level lower than $g_{j}^{lt}\left(a\right)$ leads to a non-similarity situation.  
\item An identical procedure to the one presented in the previous step can lead to obtain $g_{j}^{lt^{\prime}}\left(a\right)$. In this case, consider two performance levels, $g_{j}^{l}(b)$ and $g_{j}(a)$, with $g_{j}\left(a\right)>g_{j}^{l}(b)$. Then, in 2(a), the analyst must incrementally increase $g_{j}(a)$, and in 2(b), the
analyst must incrementally decrease $g_{j}(a)$. Figure \ref{Fig:Thresholds_1} shows the two co-constructed performance levels, $g_{j}^{lt}\left(a\right)$
and $g_{j}^{lt^{\prime}}\left(a\right)$.
\label{enu:An-identical-procedure}

\begin{figure}[htb!]
\centering 
\includegraphics{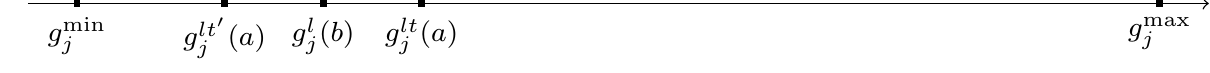} 
\caption{Co-constructed performance levels $g_{j}^{lt}\left(a\right)$ and $g_{j}^{lt^{\prime}}\left(a\right)$.}
\label{Fig:Thresholds_1}
\end{figure}

\item Consider a performance level in the upper part of the criterion scale,
say $g_{j}^{u}\left(b\right)$. A procedure as described in Steps
\ref{enu:Consider-two-performance} and \ref{enu:An-identical-procedure}
can be applied. This performance level, $g_{j}^{u}\left(b\right)$,
must be significantly different from $g_{j}^{l}(b)$. In general,
a \textquotedblleft good\textquotedblright{} choice for such a level
is in between $2/3$ to $3/4$ of the criterion scale. Figure \ref{Fig:Thresholds_2} shows the fixed levels until this step. \label{enu:Consider-a-performance}
 
\begin{figure}[htb!]
\centering 
\includegraphics{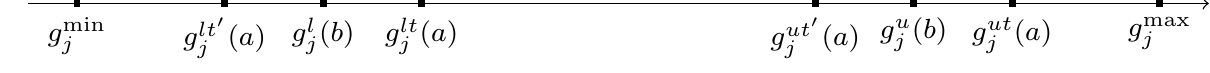} 
\caption{Co-constructed performance levels}
\label{Fig:Thresholds_2}
\end{figure}

\item Analogously, an identical procedure as previously described in Steps \ref{enu:Consider-a-criterion}-\ref{enu:Consider-a-performance}
can be applied with the aim of determining how many (if any) levels
of action $a$ must be removed from (and added to) the reference scale level to consider that both actions, $a$ and $b$, are dissimilar according to criterion $g_{j}$. The performance level $g_{j}^{ut^{\prime}}\left(a\right)$ is the
lowest constructed performance level allowing to conclude about the
dissimilarity with respect to $g_{j}^{u}(b)$. Any other performance
level lower than $g_{j}^{ut^{\prime}}\left(a\right)$ leads to a dissimilarity
situation. \label{enu:Analogously,-an-identical}
\item Finally, a procedure as in Step \ref{enu:Analogously,-an-identical}
must be followed to determine the performance levels of action $a$, for which actions $a$ and $b$ are considered totally dissimilar, according to criterion $g_{j}$.
\end{enumerate}

\paragraph{Computational aspects} The analyst must wonder whether the value of each threshold depends on the chosen reference scale level. For that, after application of the proposed interaction protocol, the difference between each fixed performance level and the respective reference level is determined. For instance, when considering two performances, $g_{j}^{lt}\left(a\right)$ and $g_{j}^{l}(b)$, we can define the performance difference, $\Delta_{j}^{lt}(a,b)$, where $\Delta_{j}^{lt}(a,b)=g_{j}^{lt}\left(a\right)-g_{j}^{l}(b$
for cardinal levels, and $\Delta_{j}^{lt}(a,b)$ is equal to the number of levels in between the two performance levels, $g_{j}^{lt}\left(a\right)$ and $g_{j}^{l}(b)$, for ordinal or discrete levels. When we obtain the same value for the difference of performances whether considering a reference level in the lower or upper part of the scale (e.g., $\Delta_{j}^{lt}(a,b)=\Delta_{j}^{ut}(a,b)$), we can conclude that the respective threshold is constant along the scale. For example, in the case of the (non-negative) similarity threshold $t(g_{j}\left(b\right))$, such a value can be determined as follows: $t\left(g_{j}\left(b\right)\right)=|\Delta_{j}^{lt}(a,b)|=|\Delta_{j}^{ut}(a,b)|$. In an identical way, the values of the other thresholds can be determined. Otherwise, an affine function, $g_{j}\left(b\right)\alpha_{j}+\beta_{j}$,
can be used to represent the variation of the threshold
along the scale. The coefficients of the affine function can be computed from the following system:

\[
\left\{ \begin{array}{l}
|\Delta_{j}^{lt}(a,b)|=g_{j}^{l}\left(b\right)\alpha_{j}+\beta_{j}\\[1mm]
|\Delta_{j}^{ut}(a,b)|=g_{j}^{u}\left(b\right)\alpha_{j}+\beta_{j}\\[1mm]
\end{array}\right.
\]

This is also valid for the other thresholds. 

\paragraph{Illustrative example} 
Let us analyze the following simple numerical example:
\begin{enumerate}
\item Consider a criterion, $g_{j},$ to be maximized and suppose that its
scale, $E_{j}$, is quantitative and continuous, with $E_{j}=[0,200].$
Let $b$ represent a reference action, and $g_{j}(b)=70$ and $g_{j}(b)=140$
be the performances chosen by the  as the two representative
levels on criterion $g_{j}$. Analyst ask DM: how much
levels can action $a$ be less (and greater) than action $b$ to consider
any similarity between these two actions on criterion $g_{j}$? Suppose
that the following information is given by the DM:
\begin{itemize}[label={--}]
\item For $g_{j}(b)=70$, the maximal
``negative difference'' and the maximal ``positive difference''
compatible with similarity between two actions is $10$ for both. 
\item For $g_{j}(b)=135$, the maximal
``negative difference'' and the maximal ``positive difference''
compatible with similarity between two actions are $20$ and $25$,
respectively.
\end{itemize}
\item In order to determine the form of the two similarity thresholds, $t_{j}\big(g_{j}(b)\big)$
and $t_{j}^{'}\big(g_{j}(b)\big)$, and taking into consideration
the information above, the following two systems can be obtained:

\begin{equation} 
\begin{cases}
10=70\,\alpha_{t^{'}}+\beta_{t^{'}}\\
20=135\,\alpha_{t^{'}}+\beta_{t^{'}}
\end{cases}\label{eq:system_t_prime}
\end{equation}
\begin{equation} 
\begin{cases}
10=70\,\alpha_{t}+\beta_{t}\\
25=135\,\alpha_{t}+\beta_{t}
\end{cases}\label{eq:system_t}
\end{equation}
\item The solutions of Systems \ref{eq:system_t_prime} and \ref{eq:system_t}
are $\alpha_{t^{'}}=\frac{2}{13}$ and $\beta_{t^{'}}=-\frac{10}{13}$,
and $\alpha_{t}=\frac{3}{13}$ and $\beta_{t}=-\frac{80}{13}$, respectively.
Therefore, we obtain the following variable similarity thresholds:

\begin{itemize}[label={--}]
\item $t_{j}^{'}\big(g_{j}(b)\big)$ $=\frac{2}{13}$ $\,g_{j}(b)-\frac{10}{13}$.
\item $t_{j}\big(g_{j}(b)\big)=\frac{3}{13}$$\,g_{j}(b)$
$-\frac{80}{13}$.
\end{itemize}

\end{enumerate}

In the same way, the two dissimilarity thresholds, $u_{j}\big(g_{j}(b)\big)$
and $u_{j}^{'}\big(g_{j}(b)\big)$, and the two total dissimilarity
thresholds, $v_{j}\big(g_{j}(b)\big)$ and $v_{j}^{'}\big(g_{j}(b)\big)$,
can be determined. Let us only consider the case of dissimilarity
thresholds.

\begin{enumerate}
\item Analyst ask DM for determining: how much levels can action
$a$ be less (and greater) than action $b$ to consider any dissimilarity
between these two actions on criterion $g_{j}$? Suppose that the
following information is provided:
\begin{itemize}[label={--}]
\item For $g_{j}(b)=70$, the minimal
``negative difference'' and the minimal ``positive difference''
compatible with dissimilarity between two actions is $20$ for both.
\item For $g_{j}(b)=135$, the minimal
``negative difference'' and the minimal ``positive difference''
compatible with dissimilarity between two actions is $40$ for both.
\end{itemize}
\item Therefore, we obtain the following two systems:

\begin{equation}
\begin{cases}
20=70\,\alpha_{u^{'}}+\beta_{u^{'}}\\
40=135\,\alpha_{u^{'}}+\beta_{u^{'}}
\end{cases}\label{eq:system_u_prime}
\end{equation}

\begin{equation}
\begin{cases}
20=70\,\alpha_{u}+\beta_{u}\\
40=135\,\alpha_{u}+\beta_{u}
\end{cases}\label{eq:system_u}
\end{equation}
\item The solutions of Systems \ref{eq:system_u_prime}and \ref{eq:system_u}
are obviously the same: $\alpha_{u^{'}}=\alpha_{u}=\frac{4}{13}$
and $\beta_{u'}$=$\beta_{u}=-\frac{20}{13}$ . Therefore, we obtain
the following variable dissimilarity thresholds:
 $$u_{j}^{'}\big(g_{j}(b)\big)=u_{j}\big(g_{j}(b)\big)=\frac{4}{13}\,g_{j}(b)-\frac{20}{13}.$$

\end{enumerate}

This means that we have the same dissimilarity thresholds whenever
the performance of action $a$ is greater or less than the performance
of action $b$. The same procedure can be applied to determining the
variable total dissimilarity thresholds.

\subsubsection{Modeling and determining the strength of similarity-dissimilarity \label{subsec:Modeling-and-determining}}
\noindent In this elicitation part,  we focus on determining the intensities of the similarity and dissimilarity between two actions, according to a given
criterion. The preference information leading to compute such intensities is be modeled through a procedure based on a deck of cards technique as the one in \citet{figueira2002determining}. Indeed, we aim at
determining the form of each one of the four components of the SD function on criterion $g_{j}$, i.e., $f_{j}^{1}$, $f_{j}^{2}$, $f_{j}^{3}$, and $f_{j}^{4}$, which are defined in the following
intervals, respectively: 
\begin{itemize}[label={--}]
\item $f_{j}^{1}$: $]-v_{j}^{'}\left(g_{j}(b)\right),-u_{j}^{'}\left(g_{j}(b)\right)[$;
\item $f_{j}^{2}$: $]-t_{j}^{'}\left(g_{j}(b)\right),0[$;
\item $f_{j}^{3}$: $]0,t_{j}\left(g_{j}(b)\right)[$;
\item $f_{j}^{4}$: $]u_{j}\left(g_{j}(b)\right),v_{j}\left(g_{j}(b)\right)[$.
\end{itemize}

\paragraph{Gathering the preference information} For determining the strength of similarity and dissimilarity judgments on a given criterion, we propose that the analyst applies a simple interaction protocol with the DM. Let us consider a criterion $g_{j}$ to be maximized, expressed on a discrete or continuous scale. A reference scale level in the central part of the criterion scale, $g_{j}^{c}(b)$, must be chosen and, taking into account the thresholds previously defined, the four intervals are determined. For each one of these intervals, let $\Delta_{j}^{1}$ denote the value of the lower bound of the interval and $\Delta_{j}^{p}$ the value of its upper bound. For instance, for $f_{j}^{3}$, $\Delta_{j}^{1}=0$ and $\Delta_{j}^{p}=t_{j}$. Now, we must consider some values in between these two (i.e., some performance differences of scale levels and $g_{j}^{c}(b)$, allowing thus to form a sequence of ordered values representing the differences of two performance levels, $\Delta_{j}^{1}\ldots,\Delta_{j}^{k},\ldots,\Delta_{j}^{p}$. The number $p$ to be considered depends on the scale. The deck of cards technique works as follows, through interaction between the analyst and the DM:

\begin{enumerate}
\item A set of $p$ cards, with values $\Delta_{j}^{1}\ldots,\Delta_{j}^{k},\ldots,\Delta_{j}^{p}$, is provided to the DM. If the she/he considers that some (consecutive) values of the performance difference are equally similar or dissimilar, then the she/he must place the corresponding cards in the same position in the ranking, building a subset of cards. At this point, the DM
has a total order of subsets, say $S^{1},\ldots,S^{k},\ldots,S^{r}.$
\item A large enough set of blank cards is provided to the DM. The similarity-dissimilarity of two successive positions of the cards (or two subsets of cards)
in the ranking can be more or less close. These blank cards are used to model the more or less \textquotedblleft closeness\textquotedblright{} of the positions. Indeed, the DM is asked to introduce the blank cards in such a way that the greater the difference between two consecutive positions, the greater the number of blank cards: no blank card means that the difference would be minimal; one blank card means twice the minimal difference; two blank cards means three times the minimal difference, and so on.
\end{enumerate}

The proposed technique must be applied to each one of the four considered intervals. Thus, we obtain the necessary preference information to compute the values of the similarity and dissimilarity intensities. Then, we can determine the form of each one of the four components.

\paragraph{Computational aspects} 
The similarity and dissimilarity intensities take values between $-1$ and $1$. While for components $f_{j}^{2}$ and $f_{j}^{3}$, similarity intensity values are assigned to the levels, within the range $[0,1]$, for $f_{j}^{1}$ and $f_{j}^{4}$, dissimilarity intensity values are assigned, within the range $[-1,0]$. In the case of a discrete
scale, when we compute the intensities, a discrete SD function is then obtained. In the case of a continuous scale, even though various forms can be assumed by these functions, in what follows, we will only consider the linear case. For assigning the intensities using
the judgments of the DM obtained through the application of the proposed deck of cards technique, we introduce a procedure based on the one proposed in \citet{bottero2018choquet} for building interval scales. We present a step-by-step procedure to compute the intensities as follows.

\begin{enumerate}
\item Consider the two subsets $S^{1}$ and $S^{r}$, and assign intensity
values to them: 
\begin{enumerate}
\item Similarity cases ($s_{j}(a,b)$): 
\begin{enumerate}
\item $f_{j}^{2}$ : $f_{j}(S^{1})=0$ and $f_{j}(S^{r})=1$; 
\item $f_{j}^{3}$ : $f_{j}(S^{1})=1$ and $f_{j}(S^{r})=0$; 
\end{enumerate}
\item Dissimilarity cases ($d_{j}(a,b)$): 
\begin{enumerate}
\item $f_{j}^{1}$ : $f_{j}(S^{1})=-1$ and $f_{j}(S^{r})=0$; 
\item $f_{j}^{4}$ : $f_{j}(S^{1})=0$ and $f_{j}(S^{r})=-1$. 
\end{enumerate}
\end{enumerate}
\item Let $e^{k}$ denote the number between two consecutive subsets, $S^{k}$
and $S^{k+1}$, for $k=1,\ldots,r-1$. Consider the following ranking: 

\[
S^{1}e^{1}S^{2}e^{2}\cdots S^{k}e^{k}S^{k+1}\cdots S^{r-1}e^{r}S^{r}.
\]
\item Compute the unit, $\alpha$, as follows: 

\[
\alpha=\frac{1}{h},
\]
where 
\[
h=\sum_{i=1}^{r-1}(e^{i}+1). \displaystyle \](the number of units between subsets $S^{1}$ and $S^{r}$).
\item Let $\Delta_{j}^{k}$ denote a value belonging to the subset $S^{k}$,
and compute the intensity for each performance difference, $f(\Delta_{j}^{k})$,
for $k=1,\ldots,r$, as follows: 
\begin{enumerate}
\item  $f_{j}^{1}$ and $f_{j}^{2}$: 
\[
f_{j}(\Delta_{j}^{k})=f_{j}(S^{1})+\alpha\left(\sum_{i=1}^{k-1}(e^{i}+1)\right),{\displaystyle \;\mbox{for}\;k=2,\ldots,r.}
\]
\item $f_{j}^{3}$ and $f_{j}^{4}$: 
\end{enumerate}

\[
f_{j}(\Delta_{j}^{k})=f_{j}(S^{1})-\alpha\left(\sum_{i=1}^{k-1}(e^{i}+1)\right),{\displaystyle \;\mbox{for}\;k=2,\ldots,r.}
\]
\end{enumerate}

All levels belonging to a given subset, $S^{k}$, for $k=1,\ldots,r$,
will have the same intensity value. Finally, we assign the intensities
values of the considered performance differences to the function of
the respective differences between each level and $g_{j}^{c}(b)$,
i.e., $f(\Delta_{j}^{k})=f(g_{j}\left(a\right)-g_{j}^{c}(b))$, for
$k=1,\ldots,r$. For instance, for the extreme values of the interval
in which component $f_{j}^{3}$ is defined, we have: $f_{j}(0)=1$,
and $f_{j}\left(t_{j}(g_{j}^{c}\left(b\right))\right)=0$.

Therefore, we assume that the function $f_{j}$ has the same form of the obtained function when considering the particular reference level $g_{j}^{c}(b)$. Then, we generalize that to the difference between the performance of two actions, $g_{j}(a)$ and $g_{j}(b)$, and obtain the SD function, $f_{j}(\Delta_{j}(a,b))=f_{j}(g_{j}(a)-g_{j}(b))$. 

In the case of continuous scales, the proposed procedure can also be applied. The value of the intensities on a given sub-interval of performance differences, between $\Delta_{j}^{s}$ and $\Delta_{j}^{t}$, with $\Delta_{j}^{s}<\Delta_{j}^{k}<\Delta_{j}^{t}$, can be defined by linear interpolation as follows:

\[
f_{j}(\Delta_{j}^{k})=f_{j}(\Delta_{j}^{s})+\frac{\Delta_{j}^{k}-\Delta_{j}^{s}}{\Delta_{j}^{t}-\Delta_{j}^{s}}\left(f_{j}(\Delta_{j}^{t})-f_{j}(\Delta_{j}^{s})\right).
\]

Hence, a piecewise-linear interpolation function is obtained. It should be remarked that other type of interpolation can be used, when it is considered more adequate according to the preferences of the DM.

\subsection{The weights of criteria} \label{subsec:Weights}
\noindent The relative importance of a given criterion is modeled by a criterion weight. These weights can be different for each category. They can be determined by using the revised Simos procedure proposed in \citet{figueira2002determining}. It should be mentioned that, at this stage, the relative importance of a criterion must be analyzed ignoring any potential interaction between criteria that could exist. Such interaction effects will  be considered later in Subsection \ref{subsec:The-interaction-coefficients}.

There are two main steps for assigning the weights of criteria, presented in the next paragraphs. 

\paragraph{Gathering the preference information} The preference information can be obtained through the following interaction protocol: 

\begin{enumerate}
\item The analyst provides a set of $n$ cards with the identification of each criterion (e.g., name, code) on each card and some additional information when necessary (e.g., scale unit, short description of the criterion);
\item The analyst asks the DM to place the cards in a ranking, from the more important to the last important one (in the case
of a tie, the DM has to place the cards in the same position of the ranking);
\item The analyst provides to the DM an enough set of blank cards and asks her/him to insert blank cards between successive positions in the ranking previously obtained. If she/he considers that a greater importance difference exists between two consecutive positions (the more the number of cards, the more the difference of importance between the subsets of criteria);
\item Finally, the analyst asks the DM to say how many times the most important subset of criteria are more important than the least important one. This is called ratio $z$.
\end{enumerate}

The above protocol should be individually applied for each predefined category when the DM considers that the relative importance of criteria differs among categories. Then, the non-normalized weights are obtained according to the algorithm described in next paragraph.

\paragraph{Computation aspects} 
In order to determine the values of the criteria weights, compute as presented below:

\begin{enumerate}
\item  Let $S^{1},\ldots,S^{k},\ldots,S^{r}$ represent the subsets of criteria cards in the same rank position in the ranking (a subset may contain only one card), where $S^{1}$ is the rank position with the least important criterion (or criteria). Let $k(S^{1}),\ldots,k(S^{k}),\ldots,k({S^{r}})$ represent the non-normalized weights of the respective subsets. Assume that $k(S^{1})=1$;
\item Let $e^{k}$ denote the number of blank cards placed between the rank position $S^{k}$ and $S^{k+1}$. Set $e$ as follows:
\[
e=\sum_{i=1}^{r}(e^{i}+1).
\]
\item Compute each weight as follows:
\[
k(S^{k})=1+u\sum_{i=1}^{k-1}e^{i}\;\;\ \mbox{with} \;\; 
u=\frac{z-1}{e}.
\]
\end{enumerate}

All criteria in the same rank position have the same weight. Then, $k_j=k(S^{k})$ for $j=1,\ldots,n$ and $k=1,\ldots,r$. These computations must be done for each category when distinct rankings are constructed by the DM. Thus, we obtain all values of $k_j^h=k(S^{k})$, for $j=1,\ldots,n$ and $h=1,\ldots,q$.

\subsection{The interaction coefficients} \label{subsec:The-interaction-coefficients}
\noindent Three types of interaction effect in pairs of criteria can be taken into account when applying the method (see \citealp{Costaetal2018}):
\begin{enumerate}
\item \textit{Mutual-strengthening effect} modeled by a positive strengthening coefficient $k_{j\ell}^{h}$, for $h=1,...,q$ (with $k_{j\ell}^{h}=k_{\ell j}^{h}$);
\item \textit{Mutual-weakening effect} modeled by a negative weakening coefficient
$k_{j\ell}^{h}$, for $h=1,...,q$ (with $k_{j\ell}^{h}=k_{\ell j}^{h}$);
\item \textit{Antagonistic effect}  modeled by a negative
 coefficient $k_{jp}^{h}$, for $h=1,...,q$.
\end{enumerate}

It should be remarked that it is assumed that an antagonistic effect in a given pair of criteria excludes the mutual interaction effects in such a criteria pair.

\paragraph{Gathering the preference information}
Firstly, the non-normalized values of the criteria weights have to be previously determined. Secondly, the analyst has to make sure that the DM has a good understanding of the interaction effects between
criteria. Finally, an interaction procedure should be followed in order to identify the possible interaction effects in some pairs of criteria and to assign values to the interaction coefficients associated with each criteria pair (see \citealt{figueira2009interaction}, for a procedure for
\textsc{Electre} methods, and see \citealt{bottero2015dealing}, for a practical application of such a procedure). In a collaborative way, the analyst and the DM may interact as follows to determine those
interaction coefficients, while considering
individually each category:

\begin{enumerate}
\item The analyst should start by questioning the DM about the potential interaction effects between two criteria that may be considered in the decision model.
For checking possible interactions in all pairs of criteria, the following analysis may be done:
\begin{enumerate}
\item Consider criterion $g_{1}$;
\item Check whether an interaction effect between $g_{1}$ and another criterion, $g_{2},g_{3},\ldots,g_{n}$, should be taken into account, identifying the criteria pair, if any;
\item Identify the type of interaction (mutual-strengthening, mutual-weakening,
or antagonistic effect) for the pairs of criteria identified in (b);
\end{enumerate}
An identical procedure is then adopted, considering in (a) each of
the remaining criteria, $g_{2},\ldots,g_{n-1}$. It is expected that the DM identifies a small number of pairs of criteria
for which there is an interaction effect.
\item The analyst should explain to the DM how the values can be assigned to the interaction coefficients associated with the respective criteria pairs. With such a purpose, the analyst should remind the DM about the meaning of each interaction effect, for instance, by illustrating that with simple examples (see, for instance, real-world applications of \textsc{Cat-SD} in \citealt{CostaChapterART}, or examples for \textsc{Electre} methods in \citealt{figueira2009interaction}).
\item The analyst must check, for each category and each criterion, the non-negative condition (see Equation \ref{eq:non-negativity} in  \ref{sec:appendix}). If it is not fulfilled, then the above
steps must be revised with the DM, and the final values assigned to the interaction coefficients must verify the non-negative condition.
\end{enumerate}

\subsection{The likeness thresholds}
\noindent An additional preference parameter is needed: a threshold must be defined for each category as the minimum likeness degree judged necessary to say that action $a$ is alike the set of reference actions $B_{h}$, for $h=1,...,q$. $\lambda^{h}$ denotes the likeness threshold of category $C_{h}$, for $h=1,...,q$. It takes a value within the range $[0.5,1]$, and it can be viewed as majority measure.


\section{Gathering the preference information: Application of the interaction protocols}
\label{sec:interaction}
\noindent In this section, we present the preference information obtained through the intervention of the analyst and the DM applying the protocols proposed in the previous section.

\subsection{Reference profiles}
\noindent In an interaction session, the analyst asked to the DM to define reference soldiers' profiles to characterize the categories previously defined. On the basis of some related documentation, and according and the experience and knowledge of the DM regarding the performance of the Army Special Forces soldiers during the training period
and while they perform their role, the DM empirically constructed
a typical and good representative profile per category, as presented in Table \ref{Table:ReferenceActions}.

\begin{table}[!htb]
\caption{Performance of the reference Army Special Forces soldiers}
\label{Table:ReferenceActions}\smallskip{}
\centering \begin{small}\resizebox{0.85\textwidth}{!}{ %
\begin{tabular}{lllllllllll}
\hline 
Category & Reference soldier & \textbf{$PF$} & \textbf{$Intel$} & \textbf{$NR$} & \textbf{$SA$} & \textbf{$MechR$} & \textbf{$VP$} & \textbf{$PmA$} & \textbf{$Pers$} & \textbf{$Med$}\tabularnewline
\hline 
Commandos & $b_{11}$  & 17.00 & 3 & 65 & 70 & 70 & 80 & 4 & 5 & 5\tabularnewline
Paratroopers & $b_{21}$  & 14.00 & 3 & 60 & 80 & 80 & 70 & 4 & 4 & 4\tabularnewline
Special Operations & $b_{31}$  & 16.00 & 4 & 70 & 70 & 70 & 75 & 4 & 4 & 4\tabularnewline
Snipers & $b_{41}$  & 15.00 & 4 & 80 & 85 & 85 & 85 & 5 & 5 & 5\tabularnewline
\hline 
\end{tabular}} \end{small} 
\end{table}

\subsection{Per-criterion SD functions}
\noindent With respect to the construction of the SD functions, we started with the first
criterion expressed on a scale with the lowest number of possible performance
levels, that is, \textit{intelligence} ($g_{2}$), $Intel$,
assessed on a five-level qualitative scale. We had expected that this could be an easy way to the DM begins reflecting about similarities and dissimilarities in pairwise comparison of performance levels. We proceeded as follows:

\begin{enumerate}
\item The analyst started by placing a card with level $3$ (satisfactory) in the table, and asking the DM if a level \textquotedblleft good", $4$ (written in another card) has some similarity with
respect to level $3$; then, the analyst asked the same for level $5$ (excellent), getting
a positive answer for both levels;
\item The analyst also asked the same question for level $2$ (weak), and level $1$ (insufficient). While for level $2$, the DM answered that some similarity exists, but being lower when comparing to $3$ and $4$, for level $1$ there is no similarity (total dissimilarity);
\item To make sure and find out whether there was coherency in the preferences of the DM,
an analogous procedure was applied with other performance levels pairs.
\end{enumerate} 

This dialogue allowed to understand that the DM favors more the similarity
for levels above the reference level than for levels below the reference
one in pairwise comparisons. 

In a quite easy way, the DM assigned values as follows:

\begin{itemize}[label={--}]
\item $1$ to the similarity when no difference between the two levels exists;
\item $0.8$ to the positive difference of one level (e.g., $g_{2}(b)=4$ and $g_{2}(a)=5$);
\item $0.6$ to the positive difference of two levels (e.g. $g_{2}(b)=3$ and $g_{2}(a)=5$);
\item $0.4$ to the negative difference of one level (e.g., $g_{2}(b)=3$ and $g_{2}(a)=2$);
\item $-0.5$ to a positive difference of at least two levels from the reference one;
\item $-1$ to a negative difference of at least two levels from the reference one, meaning that in this case totally
dissimilarity is considered and the candidate is not suitable. 
\end{itemize}

\begin{figure}[htb!]
\centering 
\includegraphics[scale=0.9]{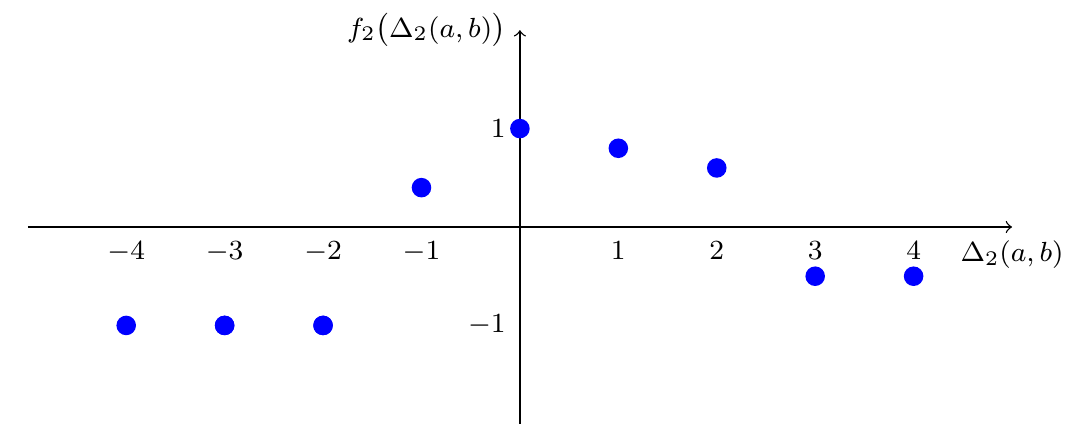} 
\caption{Per-criterion similarity-dissimilarity function for $Intel$, $PmA$, $Pers$ and $Med$}
\label{SD_Function_f2_Intel}
\end{figure}

All these preference information are represented in the piecewise function above, $f_{2}\big(\Delta_{2}(a,b)\big)$
(see also Figure \ref{SD_Function_f2_Intel}). After analyzing the cases of
$PmA$, $Pers$ and $Med$, we concluded
that the same SD function, $f_{2}$, could represent well the similarity-dissimilarity
between two candidates on these criteria (all expressed on a five-level
qualitative scale). Therefore, $f_{2}\big(\Delta_{2}(a,b)\big)=f_{7}\big(\Delta_{7}(a,b)\big)=f_{8}\big(\Delta_{8}(a,b)\big)=f_{9}\big(\Delta_{9}(a,b)\big).$

\vspace{0.25cm}
\begin{eqnarray*}
f_{2}\big(\Delta_{2}(a,b)\big)=
\begin{cases}
-1, & \mbox{\mbox{{if}}\mbox{ }}\Delta_{2}(a,b)\leqslant-2\\
\\
0.4, & \mbox{\mbox{{if}}\mbox{ }}\Delta_{2}(a,b)=-1\\
\\
1, & \mbox{\mbox{{if}}\mbox{ }}\Delta_{2}(a,b)=0\\
\\
0.8, & \mbox{\mbox{{if}}\mbox{ }}\Delta_{2}(a,b)=1\\
\\
0.6, & \mbox{\mbox{{if}}\mbox{ }}\Delta_{2}(a,b)=2\\
\\
-0.5, & \mbox{\mbox{{if}}\mbox{ }}\Delta_{2}(a,b)>2
\end{cases}
\end{eqnarray*}
\vspace{0.25cm}

As for criterion \textit{physical fitness} ($g_{1}$), $PF$, the procedure briefly described below was followed:

\begin{enumerate}
\item Starting with the level $15.00$ written in a card, the analyst asked the DM, incrementally increasing
one unit, $16.00$, $17.00$ and so on (while presenting the respective new cards), whether each level could be considered somehow similar to $15.00$. With a positive answer for $18.00,$
we got a value for a similarity threshold of $3$ ($18.00-15.00)$,
i.e., $t_{1}\big(g_{1}(b)\big)=3$;
\item Proceeding in a similar way considering
levels above $18.00$, and then below $15.00$, and generalizing for
any criteria pair with the same performance difference, the remaining
similarity-dissimilarity thresholds were assigned. 
\end{enumerate}

Then, the
intensities of the differences were determined based on rankings
of cards, as described in Subsection \ref{subsec:SD_functions}. For example, we got that a difference
of one or two levels could be considered equally similar (the cards
were placed in the same rank position). Situations of neutrality were also identified. For levels in between, it was agreed on obtained
the SD values by linear interpolation. 

The DM refereed that a difference
of six levels or more below a given reference level should be considered
as a completely dissimilarity situation, while for a difference of
five levels or more above, a quite dissimilar situation exists, but
a candidate should not necessarily be eliminated for that reason.
This was justified by the fact that the criterion is to be maximized,
and greater performances than a reference one are easily judged similar
rather than lower performances.

\begin{figure}[htb!]
\centering 
\includegraphics[scale=0.75]{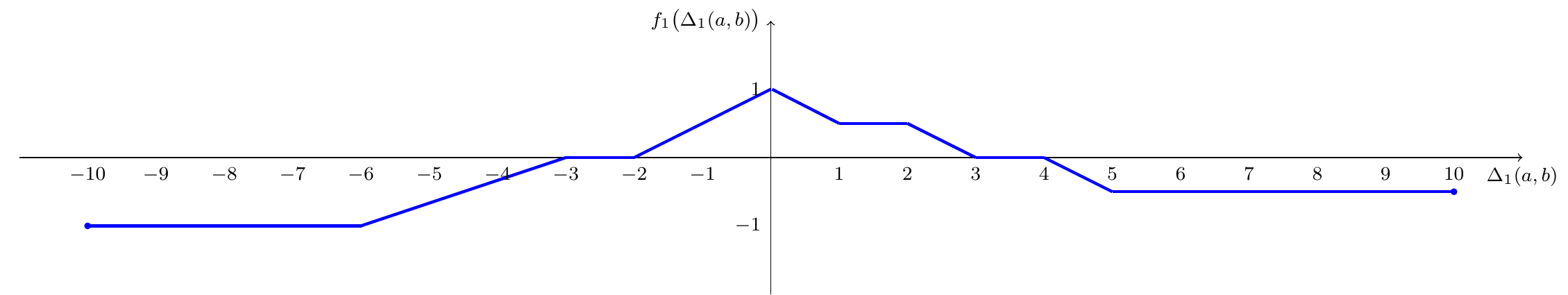} 
\caption{Per-criterion similarity-dissimilarity function for $PF$}
\label{SD_Function_f1_PF}
\end{figure}

In that way, a SD function was constructed,
as algebraically presented above by $f_{1}\big(\Delta_{1}(a,b)\big)$
and graphically represented in Figure \ref{SD_Function_f1_PF}.

\vspace{0.25cm}
\begin{eqnarray*}
f_{1}\big(\Delta_{1}(a,b)\big)=
\begin{cases}
-1, & \mbox{\mbox{{if}}\mbox{ } }\Delta_{1}(a,b)\leqslant-6\\
\\
\frac{\Delta_{1}(a,b)}{3}+1, & \mbox{\mbox{{if}}\mbox{ }}-6<\Delta_{1}(a,b)\leqslant-3\\
\\
0, & \mbox{{if}}\mbox{ }-3<\Delta_{1}(a,b))\leqslant-2\\
\\
\frac{\Delta_{1}(a,b)}{2}+1, & \mbox{{if}}\mbox{ }-2<\Delta_{1}(a,b))\leqslant0\\
\\
-\frac{\Delta_{1}(a,b)}{2}+1, & \mbox{\mbox{{if}}\mbox{ }}0<\Delta_{1}(a,b)\leqslant1\\
\\
0.5, & \mbox{\mbox{{if}}\mbox{ }}1<\Delta_{1}(a,b)\leqslant2\\
\\
-\frac{\Delta_{1}(a,b)}{2}+\frac{3}{2}, & \mbox{\mbox{{if}}\mbox{ }}2<\Delta_{1}(a,b)\leqslant3\\
\\
0, & \mbox{{if}}\mbox{ }3<\Delta_{1}(a,b))\leqslant4\\
\\
-\frac{\Delta_{1}(a,b)}{2}+2, & \mbox{{if}}\mbox{ }4<\Delta_{1}(a,b))\leqslant5\\
\\
-0.5, & \mbox{\mbox{{if}}\mbox{ } }\Delta_{1}(a,b)>5
\end{cases}
\end{eqnarray*}
\vspace{0.25cm}

As for all criteria expressed on a percentile scale ($NR$,
$SA$, $MechR$ and $VP$), it was proposed
by the DM on constructing a common function. The procedure applying the protocols described
in Subsection \ref{subsec:SD_functions} was followed.

As the level in the lowest
part of the scale the performance level $50$ was considered, and
as the level in the highest part of the scale the level $80$ was used. It was
determined that we are dealing with constant thresholds.

\begin{figure}[htb!]
\centering 
\includegraphics[scale=0.75]{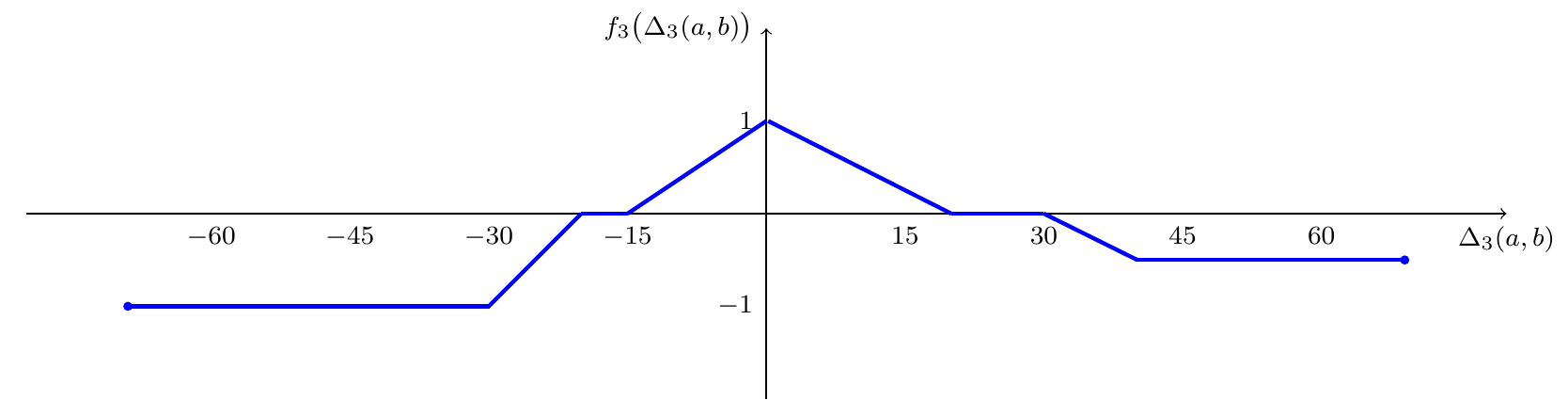} 
\caption{Per-criterion similarity-dissimilarity function for $NR$, $SA$, $MechR$ and $VP$}
\label{SD_Function_f3_NR}
\end{figure}

Having in mind the standardized percentile scales used in the psychological
tests, the DM in a quite easy way identified the SD thresholds and
the intensities of similarity-dissimilarity for performance differences.
Then, an asymmetrical SD function was built, which is algebraically
presented above by $f_{3}\big(\Delta_{3}(a,b)\big)$ and graphically
represented in Figure \ref{SD_Function_f3_NR}. It should be note that $f_{3}\big(\Delta_{3}(a,b)\big)=f_{4}\big(\Delta_{4}(a,b)\big)=f_{5}\big(\Delta_{5}(a,b)\big)=f_{6}\big(\Delta_{6}(a,b)\big)$,
in accordance with the previous statement.

\vspace{0.25cm}
\begin{eqnarray*}
f_{3}\big(\Delta_{3}(a,b)\big)=
\begin{cases}
-1, & \mbox{\mbox{{if}}\mbox{ } }\Delta_{3}(a,b)\leqslant-30\\
\\
\frac{\Delta_{3}(a,b)}{10}+2, & \mbox{\mbox{{if}}\mbox{ }}-30<\Delta_{3}(a,b)\leqslant-20\\
\\
0, & \mbox{{if}}\mbox{ }-20<\Delta_{3}(a,b))\leqslant-15\\
\\
\frac{\Delta_{3}(a,b)}{15}+1, & \mbox{{if}}\mbox{ }-15<\Delta_{3}(a,b))\leqslant0\\
\\
-\frac{\Delta_{3}(a,b)}{20}+1, & \mbox{\mbox{{if}}\mbox{ }}0<\Delta_{3}(a,b)\leqslant20\\
\\
0, & \mbox{\mbox{{if}}\mbox{ }}20<\Delta_{3}(a,b)\leqslant30\\
\\
-\frac{\Delta_{3}(a,b)}{10}+3, & \mbox{\mbox{{if}}\mbox{ }}30<\Delta_{3}(a,b)\leqslant40\\
\\
-0.5, & \mbox{\mbox{{if}}\mbox{ } }\Delta_{3}(a,b)>40
\end{cases}
\end{eqnarray*}
\vspace{0.25cm}

\subsection{Weights, interaction coefficients, and likeness thresholds}
\begin{figure}[htb!]
\centering 
\includegraphics[scale=0.4]{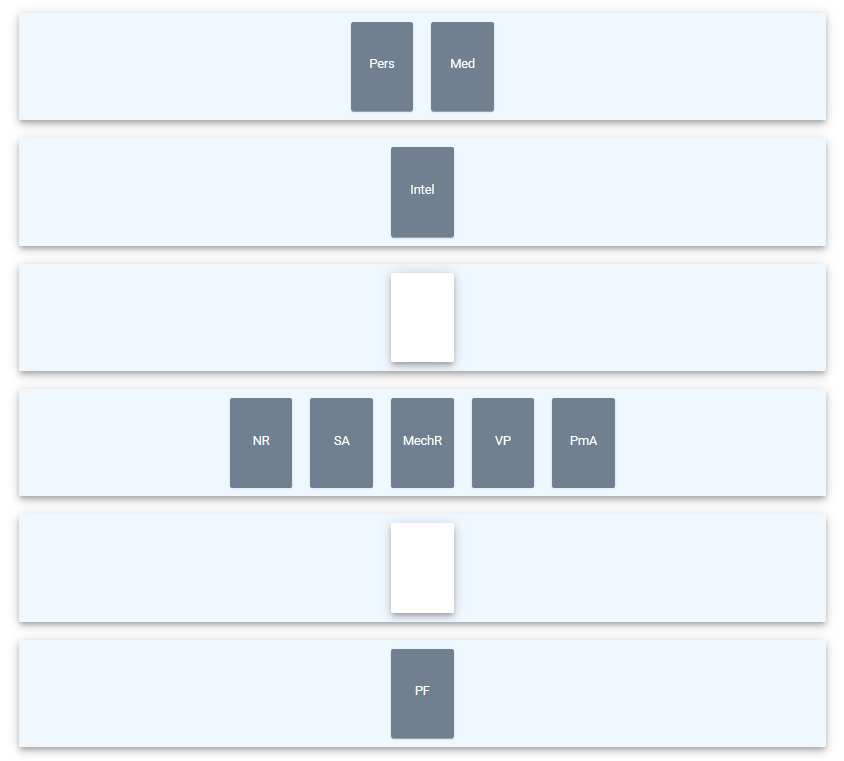} 
\caption{Cards ranking for commandos in \sc{DecSpace}}
\label{commandos_cards}
\end{figure}

\begin{figure}[htb!]
\centering 
\includegraphics[scale=0.4]{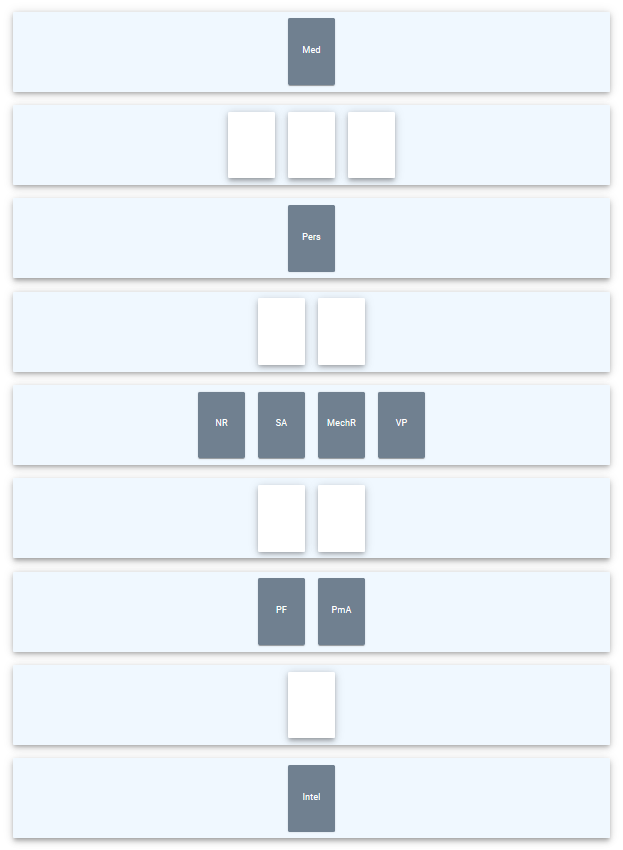} 
\caption{Cards ranking for paratroopers in \sc{DecSpace}}
\label{paratroopers_cards}
\end{figure}

\begin{figure}[htb!]
\centering 
\includegraphics[scale=0.4]{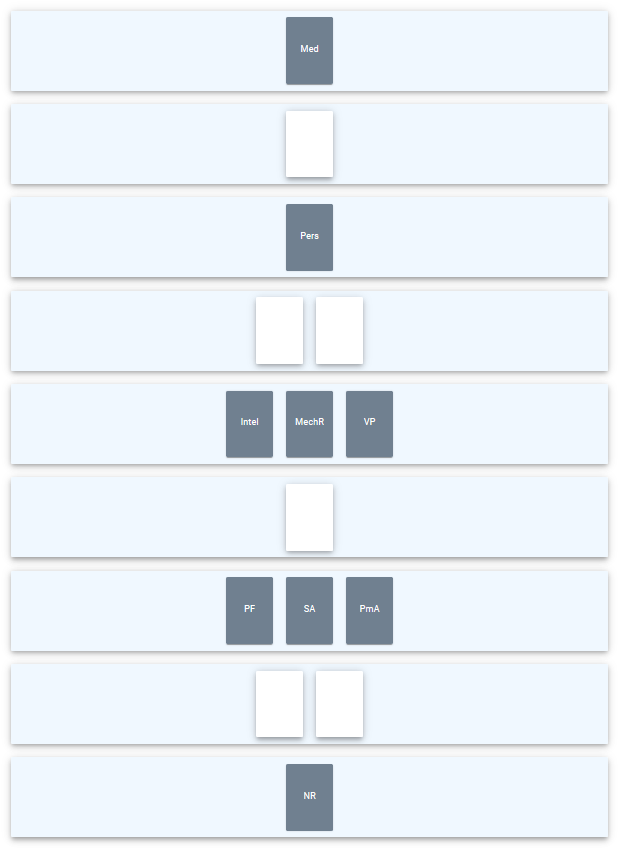} 
\caption{Cards ranking for special operations in \sc{DecSpace}}
\label{specialOperations_cards}
\end{figure}

\begin{figure}[htb!]
\centering 
\includegraphics[scale=0.4]{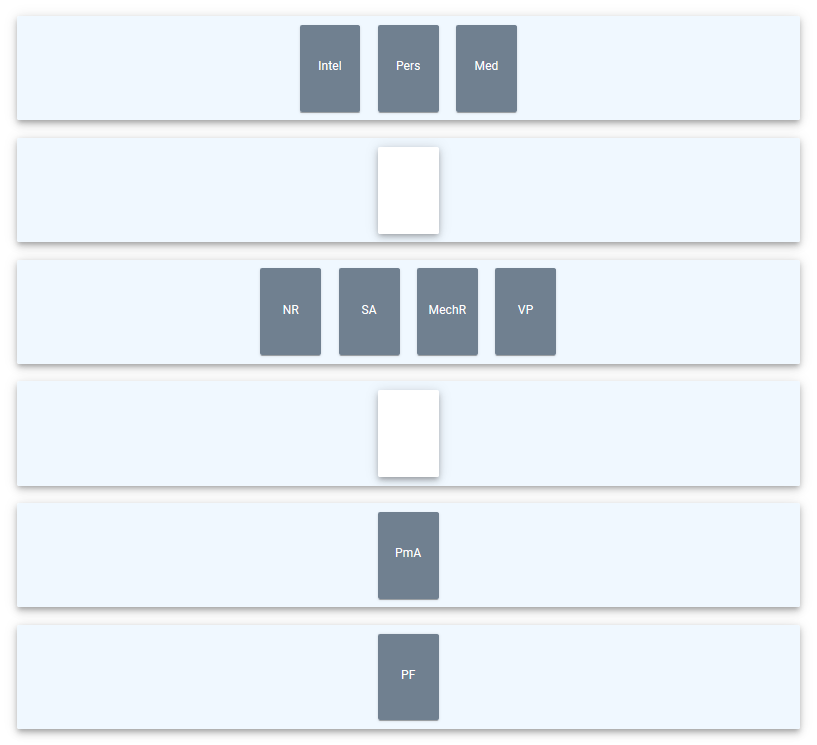} 
\caption{Cards ranking for snipers in \sc{DecSpace}}
\label{snipers_cards}
\end{figure}

\noindent During a meeting with the DM, and in the presence of a soldier of
each one of the Portuguese Army Special Forces (a commando, a paratrooper
and a special operations soldier), the protocol described in Subsection \ref{subsec:Weights}
 was applied. In a cooperative way with these three
``experts'', the DM has constructed the rankings of the cards. The
personal experience and knowledge of each soldier, specially about
the respective category, has been taken into account in this process.
Thus, a consensus of opinion among the participants allowed to obtain
the rankings and establish the ratio $z$. Some aspects and comments
expressed during the application of the protocol should be highlighted,
namely:

\begin{itemize}[label={--}]
\item \textit{General aspects}: $PF$ is a criterion with low importance
to the candidates assessment and classification, since usually in
the special training they improve their physical capabilities in a
quite easy and satisfactory way. In contrast, the criteria related
to innate abilities (e.g., $SA$, $VP$, ...) deserve greater relative
weights. The greatest importance was given to $Med$ in all categories,
justified by the fact that the soldiers will have physically demanding
training and their bodies must be ready for the challenging, without
significant medical issues; otherwise they cannot perform well;
\item \textit{Commandos} ($C_{1}$): $Med$ and $Pers$ are the
most important criteria, because the physical robustness of a commando,
and her/his emotional stability and motivation are crucial. Immediately
above is $Intel$, since planning and solving problems in unexpected
and complicated situations are important aspects to commandos. All
criteria regarding innate abilities has been considered equally
important. $PF$ is the criterion judged as the least important one.
In between consecutive levels the DM placed some blank cards and thus
constructed a final cards ranking, as depicted in Figure \ref{commandos_cards}.
Since all criteria have a similar relative importance, a rather
low value was assigned to the ratio between the weight of $Med$ and
$Pers$ (the most important criteria) and the weight of $PF$
(the least important criterion): $z=4$;
\item \textit{Paratroopers} ($C_{2}$): $Med$ is more important than the
remaining criteria. With a medical problem they cannot do the parachute
jumping. $Pers$ is the second more important. A group of
criteria has been considered equally important to paratroopers. They
are related to personal abilities: $NR$, $SA$,
$MechR$ and $VP$. Drawing logical conclusions based on numerical
data, having a great spatial ability, anticipating movements, and
acting with a great speed are relevant capabilities to successfully
jump with a parachute. Then, $PF$ and $PmA$ appear
with the same importance. $Intel$ is in the last position
of the ranking. The DM placed a certain set of blank cards between
consecutive levels, the ranking was defined, as in Figure \ref{paratroopers_cards}.
As for the value to the ratio $z$, there was a hesitation between
$6$ and $8$. They argued that there is clearly a significant difference
between the weight of most important criterion and the least one.
For this reason, initially the value $8$ was assigned. After some
reflection, they have considered the value $6$ more adequate: $z=6$; 
\item \textit{Special Operations} ($C_{3}$): $Med$ is the most important
criterion, followed by $Pers$, which is crucial when they operate.
Resilience and adaptability, among other relevant characteristics,
should be present in these soldiers. $Intel$ $MechR$
and $VP$ are equally important: solving problems, anticipating
movements and acting with great speed are important abilities to special
operations soldiers. With lower importance, $PF$,
$SA$ and $PmA$ were placed in the ranking. The last position
was occupied by $NR$. Then, a set of blank cards were adequately
placed. The final ranking is illustrated in Figure \ref{specialOperations_cards}.
The value established to the ratio $z$ is $6$, with some initially
hesitation around the value $7$. Thus, in the end, $z=6$;
\item \textit{Snipers} ($C_{4}$): $Med$ and $Pers$ are the most
important criteria. A sniper can act alone in the ground, therefore
she/he should be apt to operate in an intelligent way. Without a good
physical condition snipers cannot operate well (e.g., they must have
a great vision and physical endurance). $NR$, $SA$,
$MechR$ and $PmA$ were considered criteria with a relative
great relevance, being given the same importance for all of them.
$VP$ was placed in the above position, and $PF$ in the
last position of the ranking. Blank cards were placed according to
the preferences of the DM, as shown in Figure \ref{snipers_cards}.
The DM mentioned that all criteria have a relative closed importance
to a sniper. So, the value established to the ratio $z$ is $5$,
with some hesitation in choosing the value $4$ instead. Then, $z=5$.
\end{itemize}

Firstly, we used the DCM-SRF method available in \textsc{DecSpace}\footnote{\textsc{DecSpace} Pre-Alpha is available at http://app.decspacedev.sysresearch.org} (currently in the pre-alpha release) to obtain the values of the criteria weights, taking into
account the various values of $z$ considered by the DM. The Deck
Cards Method (DCM) is one of the early current methods supported by
\textsc{DecSpace}, mainly designed for determining the criteria weights
using the revised Simos\textquoteright{} procedure (\textquotedblleft SRF
method\textquotedblright ). Secondly, we confronted the DM with all
sets of weights to select the most adequate one for each category
(with exception of category $C_{1}$, for which a value of $z$ was
well defined, and consequently required only a final validation of
the obtained values). Finally, after the validation of the DM, the
sets presented in Table \ref{Table:Weights} were established.

\begin{table}[!htb]
\caption{Criteria weights per Special Forces category ($k_{j}^{h}$)}
\label{Table:Weights}\smallskip{}
\centering \begin{small}\resizebox{0.75\textwidth}{!}{ %
\begin{tabular}{llllllllll}
\hline 
Category & \textbf{$PF$} & \textbf{$Intel$} & \textbf{$NR$} & \textbf{$SA$} & \textbf{$MechR$} & \textbf{$VP$} & \textbf{$PmA$} & \textbf{$Pers$} & \textbf{$Med$}\tabularnewline
\hline 
Commandos  & 1 & 3.4 & 2.2 & 2.2 & 2.2 & 2.2 & 2.2 & 4 & 4\tabularnewline
Paratroopers & 1.83 & 1 & 3.08 & 3.08 & 3.08 & 3.08 & 1.83 & 4.33 & 6\tabularnewline
Special Operations  & 2.5 & 3.5 & 1 & 2.5 & 3.5 & 3.5 & 2.5 & 5 & 6\tabularnewline
Snipers & 1 & 5 & 3.4 & 3.4 & 3.4 & 1.8 & 3.4 & 5 & 5\tabularnewline
\hline 
\end{tabular}} \end{small} 
\end{table}

In order to know whether interaction in pairs of criteria could be
considered in the model, the analyst began by giving a brief explanation
to the DM about the main idea of the possible interaction effects
between two criteria. The analyst made use of examples of interactions
in a scenario related to the evaluation of cars, as presented in \Citet{corrente2016multiple}. Once a clear understanding was expressed by the DM,
through a systematic analysis, as described in Subection \ref{subsec:The-interaction-coefficients}, the following
possible interactions were identified, being considered valid for
all categories:

\begin{itemize} [label={--}]
\item[$a)$]  Mutual-strengthening effects:
\begin{itemize}
\item $PF$ ($g_{1}$) and $PmA$ ($g_{7}$)
\item $PF$ ($g_{1}$) and $Med$ ($g_{9}$)
\end{itemize}
\item[$b)$]  Mutual-weakening effects:
\begin{itemize}
\item $MechR$ ($g_{5}$) and $PmA$ ($g_{7}$)
\item $VP$ ($g_{6}$) and $PmA$ ($g_{7}$)
\end{itemize}
\end{itemize}

Making sure that the DM clearly understood the idea behind the interaction
coefficients, values were assigned to the four cases of interaction.
For the category \textit{commandos}, the following reasoning was considered to
define the final values:

\begin{itemize} [label={--}]
\item Taking into account the sum of the weights of the criteria $PF$ ($g_{1}$)
and $PmA$ ($g_{7}$), $k_{1}^{1}+k_{7}^{1}=1+2.2=3.2$, the DM assigned
a weight of $4.2$ to the coalition of $g_{1}$ and $g_{7}$. Then,
the value of this strengthening coefficient is 1 ($k_{17}^{1}=4.2-3.2=1$);
\item Taking into account the sum of the weights of the criteria $PF$ ($g_{1}$)
and $Med$ ($g_{9}$), $k_{1}^{1}+k_{9}^{1}=1+4=5$, the DM assigned
a weight of 7 to the coalition of $g_{1}$ and $g_{9}$. Then, the
value of this strengthening coefficient is $2$ ($k_{19}^{1}=7-5=2$);
\item Taking into account the sum of the weights of the criteria $MechR$
($g_{5}$) and $PmA$ ($g_{7}$), $k_{5}^{1}+k_{7}^{1}=2.2+2.2=4.4$,
the DM assigned a weight of $3.4$ to the coalition of $g_{5}$ and
$g_{7}$. Then, the value of this weakening coefficient is $-1$ ($k_{57}^{1}=3.4-4.4=-1$);
\item Taking into account the sum of the weights of the criteria $VP$ ($g_{6}$)
and $PmA$ ($g_{7}$), $k_{6}^{1}+k_{7}^{1}=2.2+2.2=4.4$, the DM
assigned a weight of $4$ to the coalition of $g_{6}$ and $g_{7}$.
Then, the value of this weakening coefficient is $-0.4$ ($k_{67}^{1}=4-4.4=0.4$).
\end{itemize}

The DM argued that the same values of the interaction coefficients
could be considered for the remaining three categories. Initially, the non-negativity
condition (see Appendix, Equation \ref{eq:non-negativity}) was not fulfilled for criterion $PmA$ ($g_{7}$), since
the initial values assigned to the weakening coefficients were very
low. Then, those coefficients were revised for all categories, and
finally they were defined as presented above. Accordingly, we have $k_{17}^{h}=1$, $k_{19}^{h}=2$, $k_{57}^{h}=-1$, and $k_{67}^{h}=-0.4$, for $h=1,...,4$.

Afterward, the analyst asked the DM for considering the assessment
of a given candidate as an overall likeness degree with respect to
a reference profile with a value between $0$ and $1$ (maximum likeness
degree). Thus, while considering individually each category, we asked
the DM how much should be such a degree to consider that the candidate
is suitable to the category under analysis. We mentioned that this
value defines a likeness threshold within the range $[0.5,1]$, being
considered as majority of votes in favor of likeness of the candidate
and the reference soldier profile. The DM refereed that: (i) for paratroopers
the minimum likeness degree is enough; (ii) for commandos and special operations
there is a relative demand on having similar soldiers to the reference
profile of the respective category (the candidates should present the required overall performance in order to be prepared for the initial intensive training); and (iii) for snipers there is a greater
requirement in terms of corresponding to the desired profile, guaranteeing that
they can fulfill their function. Accordingly, the DM defined the values
of the likeness thresholds as follows:

\begin{itemize}
\item[$-$] \textit{Commandos} ($C_1$): $\lambda^{1}=0.65$;
\item[$-$] \textit{Paratroopers} ($C_2$): $\lambda^{2}=0.50$;
\item[$-$] \textit{Special Operations} ($C_3$): $\lambda^{3}=0.65$;
\item[$-$] \textit{Snipers} ($C_4$): $\lambda^{4}=0.80$.
\end{itemize}

\section{The {\sc{Cat-SD}} in the {\sc{DecSpace}}} \label{sec:DecSpace}
\noindent In this section, we describe the design and use of \textsc{Cat-SD} in the \textsc{DecSpace}. We present first the main features and functionalities of \textsc{DecSpace}. Then, we illustrate how to make use of the available implementation of \textsc{Cat-SD} while presenting the model constructed in the case study introduced in Section \ref{sec:CAT-SD}.

\subsection{Overview of {\sc{DecSpace}}}
\noindent 
\textsc{DecSpace} is an innovative web-based platform to explore MCDA methods, conceived to be user-friendly and to offer an intuitive graphical user interface (in any web browser), designed while a certain expertise in MCDA from the user is not necessarily required. It was designed to make available several methods and to make it possible to add new ones without much programming effort. This platform is intended
for use in teaching and researching on MCDA methods, as well as for professional use as a decision support systems (DSS) for supporting decisions using those methods. \textsc{DecSpace} is currently in Pre-Alpha, meaning that the platform has been in development, and this is an early release. 

In terms of infrastructure, \textsc{DecSpace} consists of the web browser utilized by the user (where the client tier is deployed);
a dedicated server that executes most of the necessary computations;
and an external database that keeps all the data secure. Accordingly,
the platform consists of a three-tiers architecture: 

\begin{enumerate}
\item Client tier: It implements the user interface and sends user requests
to the application tier; 
\item Application tier: It confines most of the platform\textquoteright s
complexity, i.e., where most of the computational work is performed,
and carries out connections with the other two tiers; 
\item Data tier: It stores and retrieves all the data and replies to any
data requests sent by the application tier. 
\end{enumerate}

\begin{figure}[htb!]
\centering 
\includegraphics[scale=0.4]{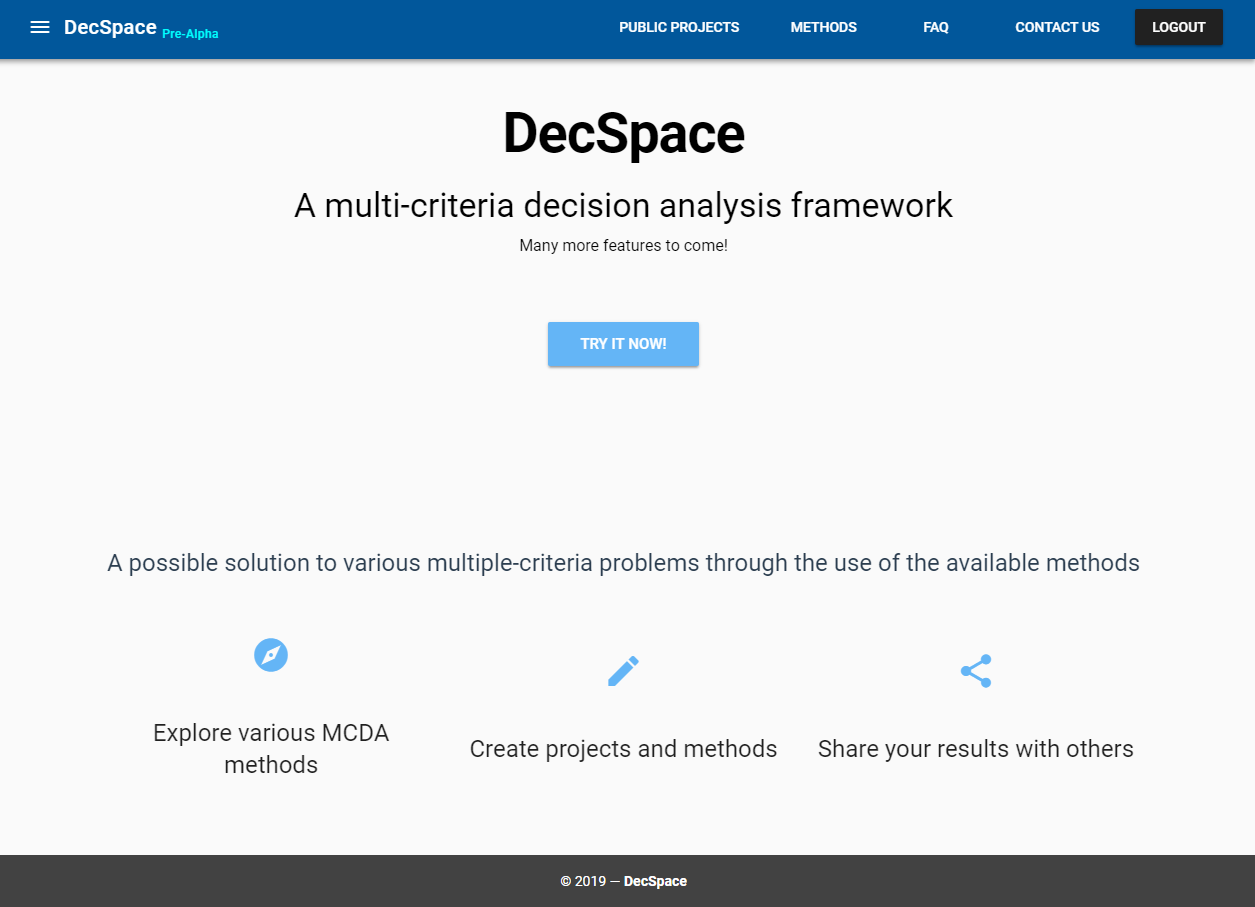} 
\caption{\textsc{DecSpace} homepage (alpha version)}
\label{Fig:decspace_homepage}
\end{figure}

A relevant inherent characteristic of \textsc{DecSpace} is that it
is a web platform accessible from devices with an Internet connection
and a browser. Its interface is thus optimized for different types
of devices, including mobile devices as tablets and smartphones. Figure \ref{Fig:decspace_homepage} displays the \textsc{DecSpace} homepage. 

To start using \textsc{DecSpace} with all the features available, a registration
process is mandatory. It is also possible to use it without
registration, as an anonymous user, but with restrictions. Any registered
user owns a project area. Each project has its own information, including the privacy
setting, public or private. The public project area provides access to projects
that where shared by other users, which can be opened, but
modifications are only allowed when the project is duplicated as a private project. Private
projects are solely available for the own user. 

\begin{figure}[htb!]
\centering 
\includegraphics[scale=0.4]{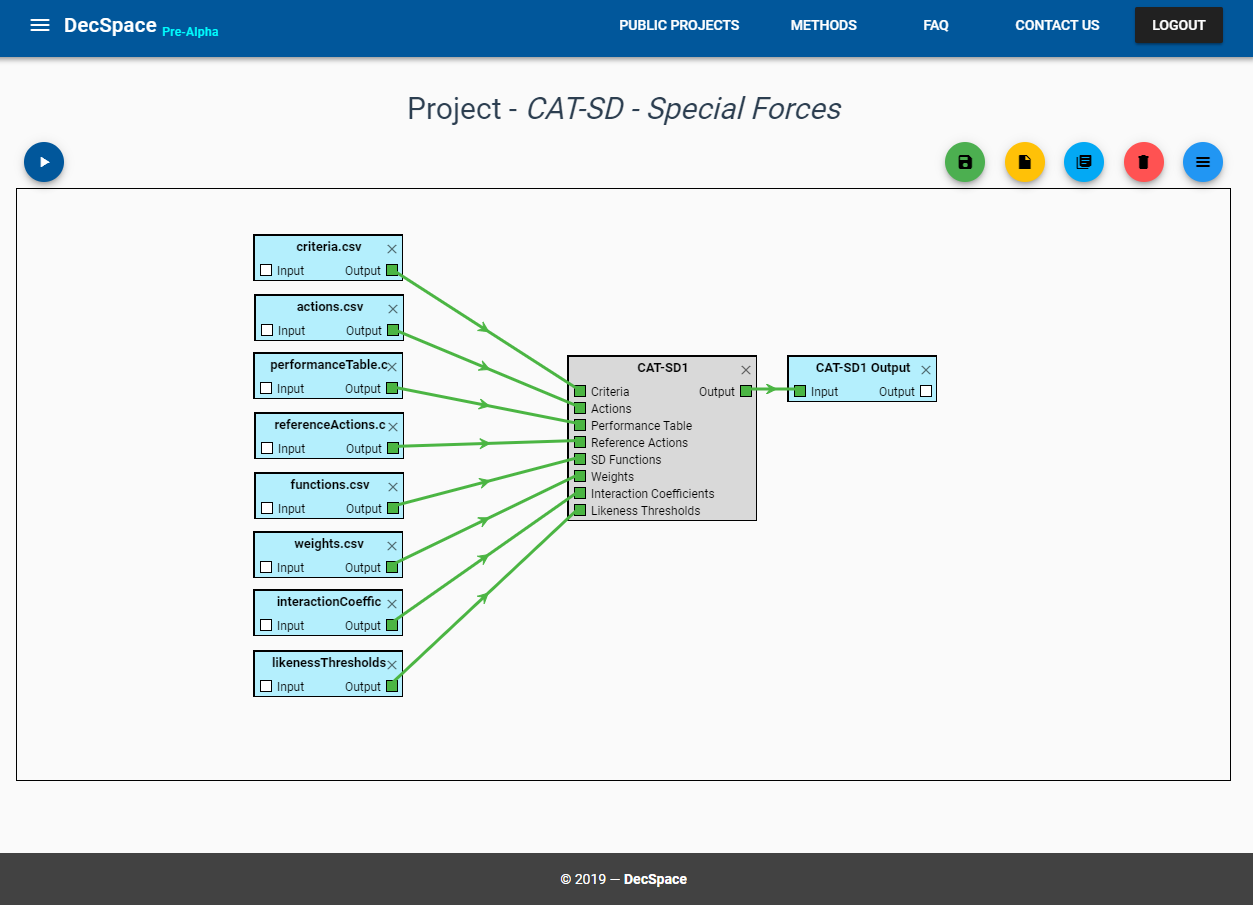} 
\caption{A \textsc{Cat-SD} workflow in \textsc{DecSpace}} \label{Fig:CAT-SD_workflow}
\end{figure}

The workspace area, as presented in Figure \ref{Fig:CAT-SD_workflow}, is intend to support the
building and running of workflows. THe MCDA methods modules can be chosen from the available methods list, and by doing that, the corresponding ``boxes''
show up in the workspace. The methods can be locally implemented by
developers, or be available as remote methods, as for example from the \textit{diviz} server \citep{MeyerBigaret2012diviz}. The user is able to manually enter
data into the MCDA methods modules or to import data files (data modules) and connect them to the method modules. Moreover, the input data can
be imported in the format of Comma-Separated Values (CSV) and JavaScript
Object Notation (JSON). Importing a .zip with several files and a .zip
with other workflows is also possible. In order to get
the results of the method modules, the workflow must be executed.

\textsc{DecSpace} currently presents the following main features and
related functionalities:

\begin{itemize}[label={--}]
\item Types of Users: There are four types of users with different permissions: 
\begin{enumerate}
\item Developers: They can implement and add new MCDA methods;
\item Administrator: She/he manages all users and projects, having permission
to modify or delete any object;
\item Anonymous users: They can test the platform and explore it with the
restriction of having temporary projects (they cannot be saved);
\item Registered users: They have an account and their work is persistent, that is, the users can create their own projects, having access to all the features that \textsc{DecSpace} offers. 
\end{enumerate}
\item User registration and login: To use the platform a user can register, providing an e-mail address and a password (a username is created based on the user' e-mail), or login to access the her/his personal area, if she/he is already registered. It is possible to enter as a guest, by clicking on ``Try it now!''; 
\item My Projects: In this area, a registered user can create projects (with a click
on the ``New Project'' button), choose whether the project is public
or private, and manage her/his projects, using a set of functionalities.
Each created project has associated the following: Name, Privacy,
Last Update, Created and Actions. Besides the functionality of creating
new projects, the following actions are available for each project:
Open Project (go to the ``Workspace''), Duplicate Project (create
a project equal to the selected project), and Delete Project (permanently
eliminate the project, after confirmation by the user). In addition,
the default order of visualizing the list of projects can be changed
by inverting the alphabetical order of the project name, and the number
of rows per page can be chosen among the four options (5, 10, 25, All);
\item Settings: In this area, it is possible to change the username and the password (exclusively available in the personal area of registered users);
\item Public Projects: This area contains all public projects that are shared
by users (available to everyone), with the following information for
each project: Name, Owner, Last Update and Created. The same possible
Actions as in ``My Projects'' are available (open, duplicate and
delete), as well as the possibility of changing the project order
by inverting the alphabetical order by ``Name'' or ``Owner'',
and choosing the number of rows per page; 
\item Methods: All the available methods are in this area. Each method has a short description, an example and a step-by-step explanation, intended to provide some helpful information to the user;
\item FAQ: It contains commonly asked questions and the respective answers related to some features of the platform; 
\item Workspace: This is the area where the users construct workflows by dragging and dropping, and properly connecting MCDA methods modules and data modules (``boxes") in an intuitive graphical user interface. The following actions are available: Execute Workflow, Save Project, Import Data, Method Selection, Delete Workflow and Project Menu (see Figure \ref{Fig:CAT-SD_workflow}). The data and the preference information can be manually provided in the methods modules. It is possible to import and export CSV and JSON files. A .zip file can also be imported, containing a workflow that was already used. The workflows can be executed, saved and deleted. 
\end{itemize}

\subsection{Building a {\sc{Cat-SD}} workflow}
\noindent In order to use \textsc{Cat-SD} in \textsc{DecSpace}, the user must
previously create a project. When opening such a project, it is displayed
the workspace, where a decision model can be constructed as a workflow.
The ``CAT-SD'' must be selected among the current available
methods (``Method Selection''), clicking on ``Add method'', and
immediately a method module, ``CAT-SD1'', appears in the workspace
area, as depicted in Figure \ref{Fig:CAT-SD_workflow}. All required data and preference parameters
need to be linked to such a method module. This can be done by manually
inserting the data, after clicking on the boxes modules, or by uploading
correctly structured CSV or JSON files, which appear as data modules,
to be then adequately connected to the method module. Figure \ref{Fig:CAT-SD_workflow} shows
all modules connections. At any time, data can be changed in the workspace area (changed data files can be exported later on).

As presented in Subsection \ref{subsec:flowchart} (see Figure \ref{Fig_flowchart}), the input include the main
data (criteria, actions and performance table) and the preference
parameters (SD functions, reference actions, weights, interaction
coefficients and likeness thresholds), which are clearly separated
in the method module. Only the interaction coefficients are not necessarily
required, since in some cases they are not present. The remaining
data are mandatory to be possible to execute the constructed
workflow and then obtain the results (method output), which appear as a data module in the workspace. By clicking on such a module (``box''),
the corresponding results can be viewed and analyzed by the user.
The tables presented in the separated spreadsheets are interconnect in such a way that when input data are provided (e.g., criteria), the necessary related information appear in the tables that depend on such data. 

In detail, the user can provide the following information:

\begin{itemize}
\item[$a)$] Data
\begin{itemize}[label={--}]
\item Criteria: In this table, seven columns appear by default:
\begin{enumerate}
\item Name: It is the name of the criterion;
\item Description: It is the respective criterion description or some related
information (it is not mandatory);
\item Direction: It corresponds to the preference direction (the user must
choose one of the two possible options: ``Maximize'' or ``Minimize'');
\item Scale Type: It is related to the kind of data of the criterion performance
levels (the user must choose ``Ordinal'' or ``Cardinal'');
\item Min: If the scale type is ``Cardinal'', then a minimum
value for the performance levels should be provided;
\item Max: If the scale type is ``Cardinal'', then a maximum
value for the performance levels should be provided;
\item Num Levels: If the scale type is ``Ordinal'', then the total
number of scale levels should be provided;
\end{enumerate}
\item Actions: Data related to the potential actions only include:
\begin{enumerate}
\item Name: It is the name of the action;
\item Description: It is the respective action description (it is not mandatory);
\end{enumerate}
\item Performance Table: The rows correspond to the actions names and the
columns corresponds to the criteria names. Performance levels on each
criterion must be provided for each action (the platform verify that
they fulfill the criteria scales characteristics); 
\end{itemize}
\item[$b)$] Preference information
\begin{itemize}[label={--}]
\item Reference Actions: For each action, the following information is needed:
the name of the category to be considered (Category); the name of
the reference action (Name) and then a performance levels have to
be fulfilled in each criterion column;
\item SD Functions: For each criterion a set of rows with a value within
the range $[-1,1]$ (SD Value) and the respective performance differences
for which the function takes the SD Value must be provided in a form
of a mathematical condition (Condition);
\item Weights: For each category and each criterion, a value of the criterion
weight must be given, i.e., a set of weights per category;
\item Interaction Coefficients: Firstly, the user can choose a category
among the predefined categories that appear as the possible options
(Category). Secondly, for such a category, a first criterion can be
chosen among the options, i.e., the previously defined criteria (Criterion
1), and then a second criterion, among the remaining criteria has
to be selected (Criterion 2). Thirdly, the type of interaction can
be chosen among the three options (Type). Finally, a value for the
interaction coefficient must be provided (Value). This procedure has
to be followed for all interaction coefficients considered in the
model. Alternatively, the data can be previously organized in a file,
and then imported and adequately connected, as for the input data
and the rest of preference parameters. The platform checks the non-negative
condition (see Equation \ref{eq:non-negativity} in \ref{sec:appendix}) and notifies the user, displaying an alert message box,
in case of the values do not fulfilled the condition;
\item Likeness Thresholds: For each predefined category (Category) a likeness
threshold must be defined with a value within the range $[0.5,1]$
(Value).
\end{itemize}
\end{itemize}

The constructed workflow can be executed when all connections are
properly done and the required data are provided (\textsc{DecSpace}
algorithms validate that and make the needed calculations for obtaining
the assignment results). The results of a \textsc{Cat-SD} workflow
are available in a data box that contains a table summarizing the
classification of the actions into the considered categories or possibly
into the ``Non-assigned'' category (it appears by default). Intermediate calculations and results are also available to the users (e.g., Maximum Likeness Degrees Per Category, Maximum Likeness Degrees Per Reference Action). The workflow, including the input files and results, can be exported (it can only be saved  by a registered user) as a .zip file containing .csv files.

Once we have designed a decision model, with all preference parameters
defined and validated by the DM, we used the data of the case study
in the construction of the corresponding \textsc{Cat-SD} workflow
in \textsc{DecSpace}. Figure \ref{Fig:CAT-SD_workflow} presents such a workflow. Data about
actions, performance table and SD functions were previously organized in spreadsheets (.csv files). Then, they were imported and properly connected to
the method module. The remaining data seemed easier and more user-friendly
to be inserted directly into tables in the workspace. Indeed, the \textsc{Cat-SD}
module has been designed with that purpose, offering a graphical user interface to provide data.

\subsection{Results and discussion}
\noindent We have obtained the assignment of the twenty
candidates (dummy, since actual data are classified) into five categories (categories representing four
Special Forces and the category \textit{unsuitable candidates}), according to
the constructed model. It is worth mentioning that the model was revised twice with the DM, and here we only present the final model and the respective associated results. Table \ref{Table:Results} shows the obtained results.

\begin{table}[!htb]
\caption{Assignment of the candidates to Special Forces}
\label{Table:Results}\smallskip{}
\centering \begin{small} \resizebox{0.9\textwidth}{!}{ %
\begin{tabular}{lccccc}
\hline 
Candidate & \textit{Commandos} & \textit{Paratroopers} & \textit{Special Operations} & \textit{Snipers} & \textit{Unsuitable Candidates}\tabularnewline
\hline 
$a_{1}$  & $\checkmark$ & $\checkmark$ & $\checkmark$ &  & \tabularnewline
$a_{2}$  & $\checkmark$ & $\checkmark$ & $\checkmark$ & $\checkmark$ & \tabularnewline
$a_{3}$  &  & $\checkmark$ & $\checkmark$ &  & \tabularnewline
$a_{4}$  &  & $\checkmark$ &  &  & \tabularnewline
$a_{5}$  &  & $\checkmark$ &  &  & \tabularnewline
$a_{6}$  & $\checkmark$ & $\checkmark$ & $\checkmark$ &  & \tabularnewline
$a_{7}$  &  &  &  &  & $\checkmark$ \tabularnewline
$a_{8}$  & $\checkmark$ & $\checkmark$ & $\checkmark$ & $\checkmark$ & \tabularnewline
$a_{9}$  & $\checkmark$ & $\checkmark$ & $\checkmark$ &  & \tabularnewline
$a_{10}$  &  &  &  &  & $\checkmark$ \tabularnewline
$a_{11}$  & $\checkmark$ &  & $\checkmark$ &  & \tabularnewline
$a_{12}$  &  & $\checkmark$ &  &  & \tabularnewline
$a_{13}$  &  & $\checkmark$ &  &  & \tabularnewline
$a_{14}$  &  &  & $\checkmark$ &  & \tabularnewline
$a_{15}$  & $\checkmark$ &  & $\checkmark$ &  & \tabularnewline
$a_{16}$  & $\checkmark$ & $\checkmark$ & $\checkmark$ &  & \tabularnewline
$a_{17}$  &  & $\checkmark$  &  &  & \tabularnewline
$a_{18}$  &  & $\checkmark$ &  &  & \tabularnewline
$a_{19}$  &  & $\checkmark$  &  &  &  \tabularnewline
$a_{20}$  & $\checkmark$ & $\checkmark$ & $\checkmark$ & $\checkmark$ & \tabularnewline
\hline 
\end{tabular}} \end{small} 
\end{table}

We can observe there are candidates only
suitable for one specif category, others are considered apt for more than one category, and others are not suitable at all for the Special Forces (unsuitable candidates).

As one can expect, a few number of candidates (only three out of twenty) are considered suitable to be assigned to category \textit{snipers}, since it is the most demanded Special Forces category considered in the model. Moreover, these candidates are adequate for all four Special Forces categories, that is, a soldier apt as a sniper is also apt to be any of the remaining forces. This is an expected result, according to the DM. On the contrary, a large number of candidates (fifteen out of twenty) are
assigned to category \textit{paratroopers}, which is the category judged as less demanded. Consequently, candidates are easily accepted as soldiers adequate to \textit{paratroopers}.

Only two out of twenty militar candidates are non-assigned to a Special Forces category, that is, their profiles, in a general way, are not considered adequate to these forces, according to the constructed model. In particular:

\begin{itemize}[label={--}]
\item $a_7$ is excluded (non-assigned), albeit there is some likeness degree with the reference profiles. Her/his performance of $3$ in criterion $Med$ could contributed for that;
\item $a_{10}$ is excluded because she/he has a performance of $2$ in $Pers$, although she/he has a high performance on the majority of criteria. The DM refereed that this is in line with the Army requirements, i.e., a candidate with a personality assessed as a low level cannot be an Army Special Forces soldier.
\end{itemize}

It should be refereed that soldier $a_{17}$ is not nearly considered suitable to be a sniper. She/he has a likeness degree of $0.78$ with respect to $b_41$ and it is necessary a degree of $0.8$ (likeness threshold for \textit{snipers}). This results specially for the performance of $4$ (and not $5$ as the reference profile) in $Med$, criterion with the highest weight, and consequently with the highest contribution to the likeness degree.

According to the opinion of the DM, these results are satisfactory and they are coherent with the requirements for the recruitment of the Portuguese Army Special Forces soldiers. The constructed model takes into account the perspective and perception of the DM about this issue, and the obtained results are aligned with the expectations of the DM.

\section{Lessons learned from practice} \label{sec:lessons}
\noindent In this study, we adopted a decision aiding constructive approach that requires an
active intervention of the DM, and from which a set of conclusions
in terms of lessons learned can be drawn. We highlight the aspects
presented below:

\begin{itemize} 
\item[$a)$] Meetings: We argue that having several short meetings, aiming at discussing
particular points, is more beneficial than having long meetings, with extensive
interaction and discussions, which usually involve a great amount
of information and require a great cognitive effort from the DM. In
addition, we consider that the kind of communication involving interaction
with people personally has advantage over alternative ways to communicate
(i.e., no presential meetings). Accordingly, the case study was conducted
during several face-to-face meetings (interaction sessions), which
allowed to understand the perceptions and preferences of the DM with
respect to the decision problem at hand, focusing on a particular
point in each session. Indeed, it seemed to be effective and efficient,
with a high engagement, giving us the chance to easy and quickly clarify
concepts, ideas, objectives, etc.;
\item[$b)$] Communication: In general, we think that the analyst should start by introducing
some main concepts and terminology related to MCDA and \textsc{Cat-SD}
(e.g., criteria, categories, reference actions or profiles), aiming
to build a decision model. Even though, we argue that adapting to
the background of the DM, and using some terms close to her/his area
of knowledge instead of the ones used in MCDA, sometimes can help
the conversation. We have proceeded accordingly, and the analyst and the
DM talked in a quite natural and easy way;
\item[$c)$] Interaction protocols: The DM expressed that the way in which the interaction
protocols were applied was naturally accepted and understood. In addition,
the DM recognized that these interactions and discussions allowed
to reflect about relevant aspects for the assessment and selection
of Special Forces soldiers. This is a positive aspect favoring the
application of the method in the context of the case study;
\item[$d)$] Difficulties: At the end of each meeting, the analyst identified with
the DM the steps of the protocols easy to understand and apply, and
those needing a great effort to be accomplished by the DM. Moreover,
the analyst discussed about the difficulties the DM felt along the
process of eliciting the preference parameters, namely on proving
preference information and assigning numerical values. The main tasks
presenting some difficulties to the DM are discussed below, considering
the definition of the preference parameters individually:

\begin{itemize} [label={--}]
\item Reference profiles: Regarding the definition of the reference soldiers' profiles,
we had the perception that in the beginning the DM had some difficulties
providing the profiles in an empirical way, i.e., based on observation
and experience, since real data could not be used, but then the DM
readily provided the performances for all desired profiles (also somehow supported by existent documentation). It was
clear that we are dealing with high demanding soldier's profiles
and, for this reason, in general, the performances on the criteria
are high and relatively close among the categories. Still, the DM
was able to construct distinct and representative profiles for each
category. Besides taking into account well established reference actions,
encouraging the DM to reflect about those representative actions based
on the people's experience can be useful, since it expresses
somehow the reality, i.e., what is typically observed in practice;
\item SD functions: The way in which we proceed allowed to know about norms
in the context of the study and get the preferences of the DM. All
those information were taken into account and modeled through the
SD functions in an adequate way. We observe that the guidance of the
analyst during the dialogue with the DM is crucial, but giving some
freedom to the DM expresses her/his subjective judgments is also important.
Pairwise comparison of some particular performance levels can be useful
to generalize the form of the function under construction. Handling
cards and observing graphical representations of scales and examples
of functions can help the DM to build these functions;
\item Weights: The deck cards method for determining the weights was immediately
accepted and understood. The reason for that is not only the simplicity
associated with the procedure of ranking the cards, but also the fact
that the DM is familiarized with methods involving handling cards
when performing a competency profile analysis. The $z$ value required
additional explanation from the analyst: ``Please consider that the
least important criterion (or criteria) has one vote. How much votes
should have the most important criterion (or criteria) of the ranking?''.
The DM understood and assigned the $z$ values, as described in Subsection
\ref{subsec:Weights}. Thus, we recommend procedures based on the deck cards method to
this aim (not exclusively to this task);
\item Interaction coefficients: The explanation using examples proved to
be effective, since the DM had no difficulties in understanding the
three possible interaction effects that we can model with \textsc{Cat-SD}.
Thus, proving some examples (fictitious or real ones) applied in
different contexts, avoiding biases, can be an adequate way to explain
the interactions between criteria. Although in this case the task
of assigning values to the coefficients was performed without presenting
difficulties, alternative elicitation protocols can be more adequate
in other cases (e.g., using more natural language, classifying each criteria interaction as ``low'', ``medium'' and ``strong'', etc.);
\item Likeness thresholds: The DM had a well understanding about the likeness
thresholds. In this case, the values were easily established by the
DM, being clearly different among the categories. This could not happen
in other real-world cases, in which the DM could have to deeply reflect
about such thresholds, and they can be considered as having the same
value for all categories.
\end{itemize}
\end{itemize}

It worth to be mentioned that, in a general way, the DM easily understood
and accepted our approach, with a prompt attitude towards collaboration
to the co-construction of the model. This has facilitated all the interaction
between the analyst and the DM for gathering of the preference information
and eliciting the parameters.

\section{Conclusions} \label{sec:Conclusions}
\noindent This paper introduces the design and implementation of \textsc{Cat-SD,}
an MCDA method recently proposed in the literature, as well as it proposes protocols to elicit the preference parameters used in the method. This paper simultaneously presents a real-world case of application of the method, while focusing on conducting the study through the
utilization of the designed protocols, and shows the main features
of using \textsc{Cat-SD} in \textsc{DecSpace}. 

\textsc{DecSpace} is an innovative user-friendly web-based platform
to explore MCDA methods. It allows the construction of workflows,
an easy exploration and combination of various methods, while
offering a range of visualization features with cutting edge technology.
We illustrate how this platform can facilitate the use of the method
by providing computational support, as well as providing an intuitive
interface, with a good visualization of data and obtained results.

The present case study permitted us to test the support given by \textsc{Cat-SD} to a real-world decision scenario related to
the candidates selection process. Moreover, it permits somehow to
validate the application of the protocols designed to facilitate the
elicitation of the preference parameters, and draw some conclusions
in terms of lessons learned. It seems to us that this case constitutes
a good example of how the method can help to make informed decisions
in this kind of contexts. 

Besides the interest and engagement of the DM during this study,
the DM has revealed interest in future collaboration. Further research
should be done in the context of this study, having the constructed
model as basis and performing a deeper analysis of the decision problem.
This can involve, for example, a systematic study of the desired profiles
for each considered category and a robustness analysis (e.g., sensitivity
analysis to changes on some preference parameters). As for \textsc{Cat-SD}
in \textsc{DecSpace}, future work relies on improving usability and
user experience issues.

From a more generic point of view, this kind of application can be relevant not only in other contexts of recruitment process, but also in contexts facing complex decision situations involving interconnected problems, with multiple stakeholders. For instance, this approach can be applied to current sharing cities issues, namely urban planning problems, energy issues or other complex economic and environmental problems, while taking into account economic, environmental and social criteria, in order to support decision making on smart cities solutions.

\section*{Acknowledgements}
\addcontentsline{toc}{section}{\numberline{}Acknowledgements}
\noindent This work was supported by national funds through Funda\c{c}{\~a}o para a Ci{\^e}ncia e a Tecnologia (FCT) with reference UID/CEC/50021/2019. The authors gratefully acknowledge the collaboration of the Portuguese Army through the \textit{Centro de Psicologia Aplicada do Ex{\'e}rcito (CPAE)} in the case study. The authors also acknowledge the graduate students who have been contributed to the implementation of the \textsc{Cat-SD} method in \textsc{DecSpace}. Ana Sara Costa acknowledges the financial support from Universidade de Lisboa, Instituto Superior T\'{e}nico and CEG\nobreakdash-\hspace{0pt}IST (PhD Scholarship). Jos\'{e} Rui Figueira  acknowledges the support from the FCT grant SFRH/BSAB/139892/2018 under POCH Program and European Union’s Horizon 2020 research and innovation program under Grant Agreement No. 691895 SHAR-LLM (“Sharing Cities”).

\vfill\newpage

\section*{Appendix}\label{sec:appendix}
\addcontentsline{toc}{section}{\numberline{}Appendix}
\noindent



\noindent 
The  method was designed based on the concepts of similarity and dissimilarity. The assignment of an action  to a given category depends on the way it compares to the reference actions of that category. The pairwise comparison takes into account
similarity and dissimilarity (subjective) judgments of the DM, which are modeled trough a similarity-dissimilarity function on each criterion. An overall measure of the likeness between such an action and the set of reference actions is then computed. The main steps of the method are presented below.

\begin{enumerate}
\item \textit{Per-criterion similarity-dissimilarity function.} A \textit{per-criterion similarity-dissimilarity function} 
can be defined as follows, assuming that, the criterion $g_j$ is to be maximized An action $a$ and a reference action $b$ are compared in terms of performance on  $g_{j}$. As input consider the difference of performances $\Delta_{j}(a,b)=diff\left\{ g_{j}(a),g_{j}(b)\right\}$
(in cardinal scale, $diff\left\{ g_{j}(a),g_{j}(b)\right\} =g_{j}(a)-g_{j}(b),$ and in ordinal scales, $diff\left\{ g_{j}(a),g_{j}(b)\right\}$ is the number of performance levels in between $g_{j}(a)$ and $g_{j}(b)$).  A \textit{per}\emph{-criterion similarity-dissimilarity function} is a real-valued function, $f_{j}\big(\Delta_{j}(a,b)\big)$, that has as output a value within the range $[-1,1]$ and presents the
following features.The function $f_{j}$ is a non-decreasing function, if $\Delta_{j}(a,b)\in[-diff\{g_{j}^{\max},g_{j}^{\min}\},0]$, where $g_{j}^{\min}$ and $g_{j}^{\max}$ are the lowest and the greatest values of the scale of criterion, respectively; while $f_{j}$ is a non-increasing function, if $\Delta_{j}(a,b)\in[0,diff\{g_{j}^{\max},g_{j}^{\min}\}]$. It is easy to see that $f_{j}>0$ iff criterion $g_{j}$ contributes to similarity, while $f_{j}<0$ iff criterion $g_{j}$ contributes to dissimilarity. Moreover, we can consider the following: (1) A \emph{per-criterion similarity function} $s_{j}(a,b)=f_{j}\big(\Delta_{j}(a,b)\big)$,
when $f_{j}\big(\Delta_{j}(a,b)\big)>0$, and $s_{j}(a,b)$=0, otherwise; and, (2) A\emph{ per-criterion dissimilarity function} $d_{j}(a,b)=f_{j}\big(\Delta_{j}(a,b)\big)$,
when $f_{j}\big(\Delta_{j}(a,b)\big)<0$, and $d_{j}(a,b)=0$, otherwise.
\item \textit{Comprehensive similarity.} A way to compute the comprehensive similarity between two
actions, $a$ and $b$, is using the function $s^{h}(a,b)$ presented in Equation \ref{eq:s(a,b)}. It takes into account the contribution
of all criteria to the similarity (i.e., the values of the per-criterion
similarity functions, referred above),
the criteria weights and the interaction coefficients. It may also
take into account some dissimilarity values derived from antagonist
effects, if they exist:

{\footnotesize{}
\begin{equation}
s^{h}(a,b)={\displaystyle \frac{1}{K^{h}(a,b)}\left(\sum_{j\in G}k_{j}^{h}s_{j}(a,b)+\sum_{\{j,\ell\}\in M^{h}}z\big(s_{j}(a,b),s_{\ell}(a,b)\big)k_{j\ell}^{h}+\sum_{(j,p)\in O^{h}}z\big(s_{j}(a,b),|d_{p}(a,b)|\big)k_{jp}^{h}\right),}\label{eq:s(a,b)}
\end{equation}
}{\footnotesize \par}

\noindent 
where{\footnotesize{}
\[
K^{h}(a,b)=\sum_{j\in G}k_{j}^{h}+\sum_{\{j,\ell\}\in M^{h}}z\big(s_{j}(a,b),s_{\ell}(a,b)\big)k_{j\ell}^{h}+\sum_{(j,p)\in O^{h}}z\big(s_{j}(a,b),|d_{p}(a,b)|\big)k_{jp}^{h},
\]
}
{\footnotesize \par}

with $z(x,y)=xy,$ for $h=1,...,q.$. The set $M^{h}$ contains all pairs of criteria, $g_{j}$ and $g_{\ell}$, for which there is mutual-strengthening or mutual-weakening, and $O^{h}$ contains all pairs of criteria, $g_{j}$ and $g_{p}$ for which there is antagonism. The following condition must be verified.

\begin{equation}
k_{j}^{h}\;\;\;\;-\hspace{-0.5cm}\sum_{\big\{\{j,\ell\}\in M^{h}\;:\;k_{j\ell}^{h}<0\big\}}\hspace{-1cm}\vert k_{j\ell}^{h}\vert\;\;-\sum_{(j,p)\in O^{h}}\hspace{-0.25cm}\vert k_{jp}^{h}\vert\;\geqslant\;0{\displaystyle ,\;\mbox{for all}\;j\in G{\displaystyle ;h=1,...,q.}}\label{eq:non-negativity}
\end{equation}

\item 
\textit{Comprehensive dissimilarity.}  A way to measure the comprehensive dissimilarity between
two actions, $a$ and $b$, is using the function $d(a,b)$ presented in Equation \ref{eq: d(a,b)}. It takes into account the contribution of all criteria to the dissimilarity between actions $a$ and $b$
(i.e., the values of the per-criterion dissimilarity functions, referred above).

\begin{equation}
{\displaystyle d(a,b){\displaystyle =\prod_{j=1}^{n}\big(1+d_{j}(a,b)\big)-1.}}
\end{equation}
\label{eq: d(a,b)}

\item \textit{Comprehensive likeness.}  The function below assesses the overall degree to which action $a$ is alike to action $b$. 

\begin{equation}
{\displaystyle \delta(a,b)=s^{h}(a,b)\big(1+d(a,b)\big)}.
\end{equation}
\label{subsec:likeness_f}

Based on the likeness degree between an action and a set of reference
actions, and according to the likeness threshold chosen for a given
category, a $\lambda-$\textit{likeness binary relation} can be defined
as follows:

\begin{equation}
aS(\lambda^{h})B_{h}\Leftrightarrow\delta(a,B_{h})\geqslant\lambda^{h}.\label{eq:sim_bin_relation}
\end{equation}

where
\[
\delta(a,B_{h})=\underset{\ell=1,...,|B_{h}|}{\max}\left\{ \delta(a,b_{h\ell})\right\} .
\]

\item 
\textit{Likeness assignment procedure.} For a $\lambda^{h}\in[0.5,1]$, for $h=1,\ldots,q$, the assignment procedure  follows the next steps. 

\begin{itemize}
\item[$i)$] Compare $a$ with $B_{h}$, for $h=1,\ldots,q$; 
\item[$ii)$] Identify $U=\{u\;:\;aS(\lambda^{u})B_{u}\}$; 
\item[$iii)$] Assign  $a$ to  $C_{u}$, for all $u\in U$;
\item[$iv)$] If $U=\emptyset$, assign  $a$ to category $C_{q+1}$.
\end{itemize}

The method provides a set of possible categories (possibly a single one) to which an action $a$ can be assigned to (it may be the dummy category, meaning that action $a$ is not suitable to the remaining categories).

\end{enumerate}

\vfill\newpage

\addcontentsline{toc}{section}{\numberline{}References}
\bibliographystyle{model2-names}
\bibliography{CAT-SD_DecSpace_References}

\begin{thebibliography}{18}
\expandafter\ifx\csname natexlab\endcsname\relax\def\natexlab#1{#1}\fi
\expandafter\ifx\csname url\endcsname\relax
  \def\url#1{\texttt{#1}}\fi
\expandafter\ifx\csname urlprefix\endcsname\relax\def\urlprefix{URL }\fi
\providecommand{\eprint}[2][]{\url{#2}}
\providecommand{\bibinfo}[2]{#2}
\ifx\xfnm\relax \def\xfnm[#1]{\unskip,\space#1}\fi
\bibitem[{Amador et~al.(2018)Amador, Costa, Rodrigues, Figueira and
  Borbinha}]{DecSpace_CAPSI}
\bibinfo{author}{Amador, J.}, \bibinfo{author}{Costa, A.S.},
  \bibinfo{author}{Rodrigues, R.}, \bibinfo{author}{Figueira, J.R.},
  \bibinfo{author}{Borbinha, J.}, \bibinfo{year}{2018}.
\newblock \bibinfo{title}{Exploring {MCDA} methods with {\sc{decspace}}}, in:
  \bibinfo{booktitle}{Proceedings of the 18th Conference of the Portuguese
  Association for Information Systems - Industry 4.0 and Information Systems},
  \bibinfo{address}{Santar{\'e}m, Portugal}.
\bibitem[{Barbosa(2017)}]{DecSpace_Andre}
\bibinfo{author}{Barbosa, A.}, \bibinfo{year}{2017}.
\newblock \bibinfo{title}{{\sc{DecSpace}}: A Multi-Criteria Decision Analysis
  Framework}.
\newblock Master's thesis. Instituto Superior T{\'e}cnico, Universidade de
  Lisboa. \bibinfo{address}{Lisboa, Portugal}.
\bibitem[{Bigaret and Meyer(2010)}]{bigaret:hal-00926569}
\bibinfo{author}{Bigaret, S.}, \bibinfo{author}{Meyer, P.},
  \bibinfo{year}{2010}.
\newblock \bibinfo{title}{{Diviz: {A}n {MCDA} workflow design, execution and
  sharing tool}}.
\newblock \bibinfo{journal}{{Newsletter of the EURO Working Group Multicriteria
  Aid for Decisions}} \bibinfo{volume}{3}, \bibinfo{pages}{10 -- 13}.
\bibitem[{Bottero et~al.(2015)Bottero, Ferretti, Figueira, Greco and
  Roy}]{bottero2015dealing}
\bibinfo{author}{Bottero, M.}, \bibinfo{author}{Ferretti, V.},
  \bibinfo{author}{Figueira, J.R.}, \bibinfo{author}{Greco, S.},
  \bibinfo{author}{Roy, B.}, \bibinfo{year}{2015}.
\newblock \bibinfo{title}{Dealing with a multiple criteria environmental
  problem with interaction effects between criteria through an extension of the
  {E}{\sc{lectre}} {III} method}.
\newblock \bibinfo{journal}{European Journal of Operational Research}
  \bibinfo{volume}{245}, \bibinfo{pages}{837--850}.
\bibitem[{Bottero et~al.(2018)Bottero, Ferretti, Figueira, Greco and
  Roy}]{bottero2018choquet}
\bibinfo{author}{Bottero, M.}, \bibinfo{author}{Ferretti, V.},
  \bibinfo{author}{Figueira, J.R.}, \bibinfo{author}{Greco, S.},
  \bibinfo{author}{Roy, B.}, \bibinfo{year}{2018}.
\newblock \bibinfo{title}{On the choquet multiple criteria preference
  aggregation model: Theoretical and practical insights from a real-world
  application}.
\newblock \bibinfo{journal}{European Journal of Operational Research}
  \bibinfo{volume}{271}, \bibinfo{pages}{120--140}.
\bibitem[{Corrente et~al.(2016)Corrente, Greco and
  S{\l}owi{\'n}ski}]{corrente2016multiple}
\bibinfo{author}{Corrente, S.}, \bibinfo{author}{Greco, S.},
  \bibinfo{author}{S{\l}owi{\'n}ski, R.}, \bibinfo{year}{2016}.
\newblock \bibinfo{title}{Multiple criteria hierarchy process for
  {E}{\sc{lectre}} {T}{\sc{ri}} methods}.
\newblock \bibinfo{journal}{European Journal of Operational Research}
  \bibinfo{volume}{252}, \bibinfo{pages}{191--203}.
\bibitem[{Costa et~al.(2018)Costa, Figueira and Borbinha}]{Costaetal2018}
\bibinfo{author}{Costa, A.S.}, \bibinfo{author}{Figueira, J.R.},
  \bibinfo{author}{Borbinha, J.}, \bibinfo{year}{2018}.
\newblock \bibinfo{title}{A multiple criteria nominal classification method
  based on the concepts of similarity and dissimilarity}.
\newblock \bibinfo{journal}{European Journal of Operational Research}
  \bibinfo{volume}{271}, \bibinfo{pages}{193--209}.
\bibitem[{Costa et~al.(2019)Costa, Lami, Greco, Figueira and
  Borbinha}]{CostaChapterART}
\bibinfo{author}{Costa, A.S.}, \bibinfo{author}{Lami, I.M.},
  \bibinfo{author}{Greco, S.}, \bibinfo{author}{Figueira, J.R.},
  \bibinfo{author}{Borbinha, J.}, \bibinfo{year}{2019}.
\newblock \bibinfo{title}{A multiple criteria approach for defining cultural
  adaptive reuse of abandoned buildings}, in: \bibinfo{editor}{Huber,
  S.~Geiger, M.}, \bibinfo{editor}{de~Almeida, A.} (Eds.),
  \bibinfo{booktitle}{Multiple {C}riteria {D}ecision {M}aking and {A}iding -
  Cases on Decision Making Methods and Models with Computer Implementations}.
  \bibinfo{publisher}{Springer}, \bibinfo{address}{Cham, Switzerland}, pp.
  \bibinfo{pages}{193--218}.
\bibitem[{Figueira et~al.(2009)Figueira, Greco and
  Roy}]{figueira2009interaction}
\bibinfo{author}{Figueira, J.R.}, \bibinfo{author}{Greco, S.},
  \bibinfo{author}{Roy, B.}, \bibinfo{year}{2009}.
\newblock \bibinfo{title}{{E}{\sc{lectre}} methods with interaction between
  criteria: {A}n extension of the concordance index}.
\newblock \bibinfo{journal}{European Journal of Operational Research}
  \bibinfo{volume}{199}, \bibinfo{pages}{478--495}.
\bibitem[{Figueira and Roy(2002)}]{figueira2002determining}
\bibinfo{author}{Figueira, J.R.}, \bibinfo{author}{Roy, B.},
  \bibinfo{year}{2002}.
\newblock \bibinfo{title}{Determining the weights of criteria in the
  {E}{\sc{lectre}} type methods with a revised {S}imos' procedure}.
\newblock \bibinfo{journal}{European Journal of Operational Research}
  \bibinfo{volume}{139}, \bibinfo{pages}{317--326}.
\bibitem[{Meyer and Bigaret(2012)}]{MeyerBigaret2012diviz}
\bibinfo{author}{Meyer, P.}, \bibinfo{author}{Bigaret, S.},
  \bibinfo{year}{2012}.
\newblock \bibinfo{title}{Diviz: A software for modeling, processing and
  sharing algorithmic workflows in {MCDA}}.
\newblock \bibinfo{journal}{Intelligent Decision Technologies}
  \bibinfo{volume}{6}, \bibinfo{pages}{283--296}.
\bibitem[{Mustajoki and Marttunen(2017)}]{mustajoki2017comparison}
\bibinfo{author}{Mustajoki, J.}, \bibinfo{author}{Marttunen, M.},
  \bibinfo{year}{2017}.
\newblock \bibinfo{title}{Comparison of multi-criteria decision analytical
  software for supporting environmental planning processes}.
\newblock \bibinfo{journal}{Environmental Modelling \& Software}
  \bibinfo{volume}{93}, \bibinfo{pages}{78--91}.
\bibitem[{Roy(1996)}]{roy1996multicriteria}
\bibinfo{author}{Roy, B.}, \bibinfo{year}{1996}.
\newblock \bibinfo{title}{Multicriteria {M}ethodology for {D}ecision {A}iding}.
\newblock \bibinfo{publisher}{Kluwer Academic Publishers},
  \bibinfo{address}{Dordrecht, The Netherlands}.
\bibitem[{Roy and Bouyssou(1993)}]{royaide}
\bibinfo{author}{Roy, B.}, \bibinfo{author}{Bouyssou, D.},
  \bibinfo{year}{1993}.
\newblock \bibinfo{title}{Aide Multicrit{\`e}re {\`a} la D{\'e}cision:
  M{\'e}thodes et Cas.}
\newblock \bibinfo{publisher}{Economica}, \bibinfo{address}{Paris, France}.
\bibitem[{Roy et~al.(2014)Roy, Figueira and
  Almeida-Dias}]{roy2014discriminating}
\bibinfo{author}{Roy, B.}, \bibinfo{author}{Figueira, J.R.},
  \bibinfo{author}{Almeida-Dias, J.}, \bibinfo{year}{2014}.
\newblock \bibinfo{title}{Discriminating thresholds as a tool to cope with
  imperfect knowledge in multiple criteria decision aiding: Theoretical results
  and practical issues}.
\newblock \bibinfo{journal}{Omega, The International Journal of Management
  Science} \bibinfo{volume}{43}, \bibinfo{pages}{9--20}.
\bibitem[{Thokala and Madhavan(2018)}]{thokala2018stakeholder}
\bibinfo{author}{Thokala, P.}, \bibinfo{author}{Madhavan, G.},
  \bibinfo{year}{2018}.
\newblock \bibinfo{title}{Stakeholder involvement in multi-criteria decision
  analysis}.
\newblock \bibinfo{journal}{Cost Effectiveness and Resource Allocation}
  \bibinfo{volume}{16}, \bibinfo{pages}{31--33}.
\bibitem[{Voinov et~al.(2016)Voinov, Kolagani, McCall, Glynn, Kragt, Ostermann,
  Pierce and Ramu}]{voinov2016modelling}
\bibinfo{author}{Voinov, A.}, \bibinfo{author}{Kolagani, N.},
  \bibinfo{author}{McCall, M.K.}, \bibinfo{author}{Glynn, P.D.},
  \bibinfo{author}{Kragt, M.E.}, \bibinfo{author}{Ostermann, F.O.},
  \bibinfo{author}{Pierce, S.A.}, \bibinfo{author}{Ramu, P.},
  \bibinfo{year}{2016}.
\newblock \bibinfo{title}{Modelling with stakeholders--next generation}.
\newblock \bibinfo{journal}{Environmental Modelling \& Software}
  \bibinfo{volume}{77}, \bibinfo{pages}{196--220}.
\bibitem[{Weistroffer and Li(2016)}]{weistroffer2016software}
\bibinfo{author}{Weistroffer, H.R.}, \bibinfo{author}{Li, Y.},
  \bibinfo{year}{2016}.
\newblock \bibinfo{title}{Multiple criteria decision analysis software}, in:
  \bibinfo{editor}{Greco, A.}, \bibinfo{editor}{Ehrgott, M.},
  \bibinfo{editor}{Figueira, J.} (Eds.), \bibinfo{booktitle}{Multiple Criteria
  Decision Analysis: {S}ate of Art Surveys}. \bibinfo{publisher}{Springer
  Science+Business Media}, \bibinfo{address}{New York, NY}.
  volume~\bibinfo{volume}{2}, \bibinfo{edition}{2nd} edition. pp.
  \bibinfo{pages}{1305--1341}.

\end{thebibliography}







\end{document}